\documentclass{amsart}[10pt]

\usepackage{epsfig}
\usepackage{graphics}
\usepackage{amsfonts}
\usepackage{amsmath}
\usepackage{latexsym}

\newtheorem{theorem}{Theorem}
\newtheorem{lemma}{Lemma}
\newtheorem{proposition}{Proposition}

\newtheorem{remark}{Remark} 
\numberwithin{equation}{section}
\begin{document}

\newcommand\ho{H\"older }

\title[Renewal of singularity sets]{Renewal of singularity sets of
statistically self-similar measures}

\author{Julien Barral and Stephane Seuret} 
\address{INRIA Rocquencourt, Domaine de Voluceau Rocquencourt, 78153
 Le Chesnay cedex, France}

\keywords{Random measures|Large deviations|Hausdorff
dimension|Self-similarity|Fractals}

\subjclass[2000]{60G57|60F10|28A78|28A80} 

\begin{abstract}
This paper investigates new properties concerning the multifractal
structure of a class of statistically self-similar measures. These
measures include the well-known Mandelbrot multiplicative cascades,
sometimes called independent random cascades. We evaluate the scale at
which the multifractal structure of these measures becomes
discernible. The value of this scale is obtained through what we call
the growth speed in H\"older singularity sets of a Borel measure. This
growth speed yields new information on the multifractal behavior of
the rescaled copies involved in the structure of statistically
self-similar measures. Our results are useful to understand the
multifractal nature of various heterogeneous jump processes.
\end{abstract}
\maketitle 

% insert the primary Maths Subject Classification number in the first
% bracket
% and the secondary ams number(s) in the second bracket
% e.g. \ams{60E20}{49G03;49F10}

%%%%%%%%%%%%%%%%%%%%%%%%%%%%%%%%%%%%%%%%%%%%%%%%
%%%%%%%%%%%%%%%%%%%%%%%%%%%%%%%%%%%%%%%%%%%%%%%%
%%%%%%%%%%%%%%%%%%%%%%%%%%%%%%%%%%%%%%%%%%%%%%%%
%%%%%%%%%%%%%%%%%%%%%%%%%%%%%%%%%%%%%%%%%%%%%%%%

\section{Introduction}
\label{intro}

This paper investigates new properties concerning the multifractal
structure of statistically self-similar measures. The class of
measures to which our results apply includes the well-known Mandelbrot
multiplicative cascades \cite{M}, sometimes called independent random
cascades. The case of another important class, the random Gibbs
measures, is treated in \cite{BSGIBBS}.
 
Multifractal analysis is a field introduced by physicists in the
context of fully developed turbulence \cite{FRPA}. It is now widely
accepted as a pertinent tool in modeling other physical phenomena
characterized by a wild spatial (or temporal)
variability~\cite{M2,BaMu,RJLV}.  Given a positive measure $\mu$
defined on a compact subset of $\mathbb{R}^d$, performing the
multifractal analysis of $\mu$ consists in computing (or estimating)
the Hausdorff dimension $d_\mu(\alpha)$ of H\"older singularities sets
$ E_\alpha^\mu$. These sets $ E_\alpha^\mu$ are the level sets
associated with the \ho exponent
$$
\label{defhmu}
h_\mu(t) = \lim_{r\to 0^+}\frac{\log\mu (B(t,r))}{\log (r)}
$$
(whenever it is defined at $t$). Thus
\begin{equation}\label{defhmu'}
E_\alpha^\mu =\{t: h_\mu(t) =\alpha\}.
\end{equation}
 Of course, these limit behaviors are numerically unreachable, both
when simulating model measures or when processing real
data. Nevertheless, this difficulty can be circumvented since the
Hausdorff dimension $d_\mu(\alpha)$ of $E^\mu_\alpha$ can sometimes be
numerically estimated by counting at scale $2^{-j}$ the number of
boxes $B$ (in a regular fine grid) such that $\mu(B)\approx
2^{-j\alpha}$. This number can formally be defined, for any scale
$j\geq 0$, $\varepsilon>0$ and $\alpha>0$, by
$$N_j^\varepsilon(\alpha)= \#\big \{ I\in\mathcal{I}_j:
b^{-j(\alpha+\varepsilon)}\le \mu (I)\le b^{-j(\alpha-\varepsilon)}
\},
$$
where $b$ is an integer $\geq 2$ and $\mathcal{I}_j$ stands for the
set of $b$-adic cubes of generation $j$ contained in the support of
$\mu$. Then when some multifractal formalisms are fulfilled, it can be
shown that, for some $\alpha>0$, one has
\begin{equation}\label{ld}
 d_\mu(\alpha) = \lim_{\varepsilon \rightarrow 0} \lim_{j\rightarrow
+\infty} \frac{\log N_j^\varepsilon(\alpha) }{\log 2^j}.\end{equation}
This is the case for instance for the multifractal measures used as
models \cite{M}. In this frame it is thus natural to seek for
theoretical results giving estimates of the first scale from which a
substantial part of the singularity set $E_\alpha^\mu$ is discernible
when measuring the $\mu$-mass of the elements $B$ of the regular
grid. In other words, we search for the first scale $J \ge 0$ such
that for every $j\geq J$, one has $N_j^\varepsilon(\alpha) \approx
2^{jd_\mu(\alpha)}$. This is of course important for numerical
applications and modelisation.

\smallskip

The properties studied in this paper and in \cite{BSGIBBS} rely on
this question. We provide new accurate information on the fine
structure of multiplicative cascades, which bring some answers to the
above problem. This study also shows that Mandelbrot measures and
Gibbs measures have very different behaviors from the statistical
self-similarity point of view, while they cannot be distinguished by
the form of their multifractal spectra. Finally, our results are
critical tools for the Hausdorff dimension estimate of a new class of
limsup sets (see (\ref{DIM})) involved in multifractal analysis of
recent jump processes \cite{MOIBARRALNU,NOTE,NUGEN,LEVY}.

%%%%%%%%%%%%%%%%%%%%%%%%%%%%%%%%%
\subsection*{A definition of statistical self-similarity}

Let us now specify what we mean by statistically self-similar measure
in the sequel. Our point of view takes into account a structure which
often arises in the construction of random measures generated by
multiplicative processes.

Let $\mathcal{I}$ be the set of closed $b$-adic sub-hypercubes of
$[0,1]^d$. A random measure $\mu(\omega)$ on $[0,1]^d$ ($d\ge 1$) is
said to be statistically self-similar if there exist an integer $b\ge
2$, a sequence $Q_n(t,\omega)$ of random non-negative functions, and a
sequence of random measures $(\mu^I)_{I\in \mathcal{I}}$ on $[0,1]^d$
such that:

\smallskip
\noindent
1.  For every $I\in \mathcal{I}$ and $g\in C([0,1]^d,\mathbb{R})$ one
 has ($\stackrel{d}{\equiv}$ means equality in distribution)
$$ 
\int_{[0,1]^d}g(u)\mu(\omega)(du)\stackrel{d}{\equiv}
\int_{[0,1]^d}g(u)\mu^I(\omega)(du).
$$

\smallskip
\noindent
2. With probability one, for every $I\in \mathcal{I}$ and $g\in
C(I,\mathbb{R})$ one has
\begin{equation}\label{copie}
\int_{I}g(v)\mu(\omega)(dv)= \lambda(I)\int_{I}Q_n(v,\omega)\,
g(v)\mu^I(\omega)\circ f_I^{-1}(dv),
\end{equation}
where $f_I$ stands for a similitude that maps $[0,1]^d$ onto $I$ and
$\lambda$ is the Lebesgue measure.

\smallskip
%%%%%%%%%%%%%%%%%%%%%%%%%%%%%%
Property 1. asserts that the measures $\mu^I$ and $\mu$ have the same
probability distribution. Property 2. asserts that, up to the density
$Q_n$, the behavior of the restriction of $\mu$ to $I$ is ruled by the
rescaled copy $\mu^I$ of $\mu$.

Of course, the random density $Q_n(t,\cdot)$ plays a fundamental role,
both in the construction of the measure $\mu$, which is often equal to
the almost sure weak limit of $Q_n(t,\cdot)\cdot \lambda$, and in the
local behavior of $\mu$.

We restrict ourselves to measures with support equal to $[0,1]^d$. Up
to technical refinements, our point of view can easily be extended to
measures which support is the limit set of more general iterated
random similitudes systems (\cite{HOWA,FALC,Olsen1,ArPa,B2}).

Two main classes of measures illustrate conditions 1. and 2. The first
one appears in random dynamical systems as Gibbs measures
\cite{KAKI,FAN2,BCM}. The other one consists in some
$[0,1]^d$-martingales (in the sense of \cite{K1,K2}) considered in
\cite{BM2}. This second class of measures is illustrated in particular
by independent multiplicative cascades \cite{M,KP} as well as compound
Poisson cascades \cite{BM1} and their extensions \cite{BACRYMUZY,BM2}.
As claimed above, these two classes are quite identical regarding
their multifractal structure in the sense that any measure in these
classes is ruled by the so-called multifractal formalisms
\cite{BRMICHPEY,OLSEN}.
%, which gives a description of the size of the
%H\"older singularity sets
However, the study of their self-similarity properties reveals the
notable differences that exist between them. These differences are
consequences of their construction's schemes: For random Gibbs
measures, the copies $\mu^I$ only depend on the generation of $I$
while they are all different for $[0,1]^d$-martingales (they depend on
the interval $I$).

This difference is quantitatively measured thanks to a notion which is
related with the multifractal structure, namely the {\em growth speed}
in the $\mu^I$ H\"older singularities sets $E^{\mu^I}_\alpha$ (see
Theorem A and B below). This quantity is precisely defined and studied
in the rest of the paper for independent random cascades. It yields an
estimate of the largest scale from which the observation of the
$\mu^I$'s mass distribution accurately coincides with the prediction
of the multifractal formalism.

%%%%%%%%%%%%%%%%%%%%%%%%%%%%%%%%%%%%%%%%%%%%%%%%
%%%%%%%%%%%%%%%%%%%%%%%%%%%%%%%%%%%%%%%%%%%%%%%%
%%%%%%%%%%%%%%%%%%%%%%%%%%%%%%%%%%%%%%%%%%%%%%%%
\subsection*{New limsup-sets and conditioned ubiquity}

As claimed above, the growth speed in H\"older singularities sets
naturally appears in the computation of the Hausdorff dimension a new
type of limsup-sets, which are themselves related to some
heterogeneous jump processes. Particular cases of the jump processes
considered in \cite{MOIBARRALNU,LEVY} are for instance the sum of
Dirac masses $\sum_{j\ge 0}\sum_{0\le k\le b^j-1}j^{-2}\mu
([kb^{-j},(k+1)b^{-j}])\delta_{kb^{-j}}$, and the L\'evy process $X$
in multifractal time $\mu$ defined as $\big (X\circ \mu ([0,t])\big
)_{0\le t\le 1}$.

Let $\mu$ be a statistically self-similar measure whose support is
$[0,1]$. Let $\{x_n\}$ denote a countable set of points and
$(\lambda_n)_{n\ge 1}$ be a sequence decreasing to 0 such that
$\limsup_{n\to\infty}B(x_n,\lambda_n)=[0,1]$. Typically we think as
$\{x_n\}$ as being the sequence of jump points of the processes of
\cite{MOIBARRALNU,LEVY}. It appears that the multifractal nature of
these processes is closely related to the computation of the Hausdorff
dimension of the limsup-sets $K(\alpha,\xi)$ defined for $\alpha> 0$
and $\xi>1$ and for some sequence $(\varepsilon_n)$ converging to 0~by
\begin{equation}\label{K}
K(\alpha,\xi)=\bigcap_{N\ge 1} \ \ \bigcup_{n\ge 1: \,\mu
([x_n-\lambda_n,x_n+\lambda_n])=\lambda_n^{\alpha+\varepsilon_n}} \
[x_n-\lambda_n^\xi,x_n+\lambda_n^\xi].
\end{equation}
Heuristically, $K(\alpha,\xi)$ contains the points that are infinitely
often close to a jump point $x_n$ at rate $\xi$ relatively to
$\lambda_n$, upon the condition that $\nu
([x_n-\lambda_n,x_n+\lambda_n])\sim \lambda_n^{\alpha}$. This
condition implies that $\nu$ has roughly a H\"older exponent $\alpha$
at scale $\lambda_n$. One of the main results of
\cite{MOIBARRALNU,MOIBARRALUBIQUITY} (see also \cite{NOTE}) is the
computation of the Hausdorff dimension of $ K(\alpha,\xi)$. The value
of this dimension is related to the free energy function $\tau_\mu$
considered in the multifractal formalism for measures in
\cite{HALSEY,BRMICHPEY}. For every $ q\in \mathbb{R}$ and for every
integer $j\geq 1$, let us introduce the quantities
\begin{equation}
\label{tau}
\tau_{\mu,j}(q)=-\frac{1}{j}\log_b\, \sum_{I\in\mathcal{I}_j}\mu (I)^q
\ \mbox{ and } \ \tau_{\mu}(q)=\liminf_{j\to\infty}\tau_{\mu,j}(q).
\end{equation}
  The Legendre transform $\tau_\mu^*$ of
$\tau_\mu$ at $\alpha>0$ is defined as
$\tau_\mu^*(\alpha):=\inf_{q\in\mathbb{R}}\alpha q-\tau_\mu(q)$.

Under suitable assumptions on $(\lambda_n)$, we prove in
\cite{MOIBARRALNU,MOIBARRALUBIQUITY} that, for all $\alpha$ such that
$\tau_\mu^*(\alpha)>0$ and all $\xi\ge 1$, with probability one,
($\dim$ stands for the Hausdorff dimension)
\begin{equation}\label{DIM}
\dim\, K(\alpha,\xi) = {\tau_\mu^*(\alpha)}/{\xi}.
\end{equation}

This achievement is a non-trivial generalization of what is referred
to as ``ubiquity properties'' (see \cite{DOD} and references therein)
of the resonant system (\cite{ARNOLD}) $\{(x_n,\lambda_n)\}$ . The
main difficulty here lies in the fact that $\mu$ may be a multifractal
measure and not just the uniform Lebesgue measure. Results on growth
speed in H\"older singularity set are determinant to obtain estimate
(\ref{DIM}).

%%%%%%%%%%%%%%%%%%%%%%%%%%%%%%%%%%%%%%%%%%%%%%%%
%%%%%%%%%%%%%%%%%%%%%%%%%%%%%%%%%%%%%%%%%%%%%%%%
\subsection*{Growth speed in $\mu^I$'s H\"older singularity sets.}
Let $\mu$ be a statistically self-similar positive Borel measure as
described above. As we said, multifractal analysis of $\mu$
\cite{FAN2,KI,BCM,HOWA,Mol,ArPa,B2} usually considers H\"older
singularities sets of the form (\ref{defhmu'}) and their Hausdorff
dimension $d_\mu(\alpha)$, which is a measure of their size. The
method used to compute $d_\mu(\alpha)$ is to find a random measure
$\mu_\alpha$ (of the same nature as $\mu$) such that $\mu_\alpha$ is
concentrated on $E_\alpha^\mu\cap
E_{\tau_\mu^*(\alpha)}^{\mu_\alpha}$. This measure $\mu_\alpha$ is
often referred to as an analyzing measure of $\mu$ at $\alpha$. Then,
by the Billingsley lemma (\cite{BILL} pp 136--145), one gets
$d_\mu(\alpha)=\tau_\mu^*(\alpha)$, and the multifractal formalism for
measures developed in \cite{BRMICHPEY} is said to hold for $\mu$ at
$\alpha$.  Finally, the estimate (\ref{ld}) is a direct consequence of
the multifractal formalism (\cite{R}) for the large deviation
spectrum. Thus the existence of $\mu_\alpha$ has wide consequence
regarding the possibility of measuring the mass distribution of $\mu$
at large enough scales.

\smallskip

In this paper we refine the classical approach by considering, instead
of the level sets $E_\alpha^\mu$, the finer level sets $\widetilde
E^\mu_{\alpha,p}$ and $\widetilde E^\mu_{\alpha}$ defined for a
sequence $\varepsilon_n$ going down to 0 by
\begin{eqnarray}
 \widetilde E^\mu_{\alpha,p}\!\! & = & \Big \{t\in [0,1]^d:\forall \
n\ge p,\ b^{-n(\alpha+\varepsilon_n)}\le \mu\big (I_n(t))\big )\le
b^{-n(\alpha-\varepsilon_n)}\Big\}\\
\label{wtE} \mbox{ and }\ \ 
\widetilde E^\mu_\alpha \, & = & \bigcup_{p\ge 1}E^\mu_{\alpha,p}.
\end{eqnarray}
It is possible to choose $(\varepsilon_n)_{n\ge 1}$ so that with
probability one, for all the exponents $\alpha$ such that
$\tau_\mu^*(\alpha)>0$ one has $ \mu_\alpha(\widetilde
E^\mu_\alpha)=\Vert \mu_\alpha\Vert $.

Since the sets sequence $ \{\widetilde E^\mu_{\alpha,p}\}_{p\geq 1}$
is non-decreasing and $\mu_\alpha (\widetilde
E^\mu_\alpha)=\Vert\mu_\alpha\Vert$, we can define the growth speed of
$\widetilde E^\mu_{\alpha,p}$ as the smallest value of $p$ for which
the $\mu_\alpha$-measure of $\widetilde E^\mu_{\alpha,p}$ reaches a
certain positive fraction $f\in (0,1)$ of the mass of $\mu_\alpha$,
i.e.
$$
GS(\mu,\alpha)=\inf\left \{p: \mu_\alpha (\widetilde
E^\mu_{\alpha,p})\ge f\, \Vert \mu_\alpha\Vert \right \}.
$$ For each copy $\mu^{I}$ of $\mu$, the corresponding family of
analyzing measures $\mu^{I}_\alpha$ exists and are related with
$\mu^I$ as $\mu_\alpha$ is related with $\mu$. The result we focus on
in the following is the asymptotic behavior as the generation of $I$
goes to $\infty$ of
\begin{equation}\label{GS}
GS(\mu^{I},\alpha) = \inf\left \{p: \mu^{I}_\alpha \big (\widetilde
E^{\mu^{I}}_{\alpha,p} \big ) \ge f\, \Vert \mu^{I}_\alpha\Vert \right
\}.
\end{equation}
 This number yields an estimate of the number of generations needed to
observe a substantial amount of the singularity set
$E^{\mu^I}_{\alpha}$.  Let
$$
\mathcal{N}_n(\mu^{I},\alpha)=\#\big \{J\in\mathcal{I}_n:
b^{-n(\alpha+\varepsilon_n)}\le \mu^{I} (J)\le
b^{-n(\alpha-\varepsilon_n)}\big\}.
$$
As a counterpart of controlling $GS(\mu^{I},\alpha)$, we shall also
control the smallest rank $n$ from which $
\mathcal{N}_n(\mu^{I},\alpha)$ behaves like
$b^{n\tau_\mu^*(\alpha)}$. This rank is defined by
$$
GS'(\mu^{I},\alpha)\!=\!\inf\!\big \{p: \!\forall\ n\ge p,
b^{n(\tau_\mu^*(\alpha)-\varepsilon_n)}\le
\mathcal{N}_n(\mu^{I},\alpha)\le
b^{n(\tau_\mu^*(\alpha)+\varepsilon_n)}\big \}
$$
and yields far more precise information than a result like
(\ref{ld}).

%%%%%%%%%%%%%%%%%%%%%%%%%%%%%%%%%%%%%%%%%%%%%%%%
%%%%%%%%%%%%%%%%%%%%%%%%%%%%%%%%%%%%%%%%%%%%%%%%
\subsection*{A simplified version of the main results}

In this paper, we focus on the one-dimensional case and deal with a
slight extension of the first example of $[0,1]$-martingales
introduced in \cite{M}, called independent random cascades (see
Section~\ref{IRC}). Let us give simplified versions of the main
results detailed in Section~\ref{CASCADES}. We start with a recall of
the theorem proved in ~\cite{BSGIBBS}.

\medskip

%%%%%%%%%%%%%%%%%%%%%%%%%%%%%%%%%%%%%%%%%%%%%%%%
\noindent
{\bf Theorem A.} {\em Let $\mu$ be a random Gibbs measure as in
\cite{BSGIBBS} (in this case $\mu^I=\mu^J$ if $I$ and $J$ are of the
same generation). Suppose that $\tau_\mu$ is $C^2$. Let
$\beta>0$. There exists a choice of $(\varepsilon_n)_{n\ge 1}$ such
that, with probability one, for all $\alpha>0$ such that
$\tau_\mu^*(\alpha)>0$, if $I$ is of generation $j$ large enough, then
$ GS(\mu^{I},\alpha)\le \exp\sqrt{\beta \log j}$.}
%%%%%%%%%%%%%%%%%%%%%%%%%%%%%%%%%%%%%%%%%%%%%%%%

\medskip

The fact that $ GS(\mu^{I},\alpha)$ behaves like $o(j)$ as
$j\to\infty$ is a crucial property needed to establish (\ref{DIM}) for
random Gibbs measures.

 Theorem {\bf B} shall be compared with Theorem {\bf A}. Under
suitable assumptions, we have (see Theorem~\ref{cascades_ren})

\medskip

%%%%%%%%%%%%%%%%%%%%%%%%%%%%%%%%%%%%%%%%%%%%%%%%
\noindent
{\bf Theorem B.} {\em Let $\mu$ be an independent random cascade. Let
$\eta>0$. There exists a choice of $(\varepsilon_n)_{n\ge 1}$ such
that, with probability one, for all $\alpha>0$ such that
$\tau_\mu^*(\alpha)>0$, if $I$ is of generation $j$ large enough, then
$ GS(\mu^{I},\alpha)\le j \log^{\eta} j $.}
%%%%%%%%%%%%%%%%%%%%%%%%%%%%%%%%%%%%%%%%%%%%%%%%

\medskip

Consequently we lost the uniform behavior over $\mathcal{I}_j$ of $
GS(\mu^{I},\alpha)$ like $o(j)$, which was determinant to get
(\ref{DIM}). In fact this ``worse'' behavior is not surprising, since
Theorem B controls simultaneously $b^j$ distinct measures $\mu^I$ at
each scale $j$, while Theorem A controls only one measure at each
scale. Nevertheless, this technical difficulty can be circumvented, by
using a refinement of Theorem B (see Theorem~\ref{cascades_ren2} and
\ref{analogue2}), which is enough to get (\ref{DIM}).

\medskip

%%%%%%%%%%%%%%%%%%%%%%%%%%%%%%%%%%%%%%%%%%%%%%%%
\noindent
{\bf Theorem C.} {\em Let $\mu$ be an independent random cascade. Let
$\eta>0$. There exists a choice of $(\varepsilon_n)_{n\ge 1}$ such
that for every $\alpha>0$ such that $\tau_\mu^*(\alpha)>0$, with
probability one, for $\mu$-almost every $t$, for $j$ large enough,
$GS(\mu^{I_j(t)},\alpha)\le j\log ^{-\eta}j$.}
%%%%%%%%%%%%%%%%%%%%%%%%%%%%%%%%%%%%%%%%%%%%%%%%

\medskip

The paper is organized as follows. Section~\ref{principes} gives new
definitions and establishes two general propositions useful for our
main results. In Section~\ref{CASCADES} independent random cascades
are defined in an abstract way. This makes it possible to consider
Mandelbrot measures as well as their substitute in the critical case
of degeneracy. Then the main results (Theorems~\ref{cascades},
\ref{cascades_ren} and \ref{cascades_ren2}) are stated and
proved. Theorem~\ref{cascades_dev} is a counterpart of Theorem B in
terms of $GS'(\mu^I,\alpha)$. Theorem~\ref{cvtau} deals with a problem
connected with the estimate of the growth speed in singularities sets,
namely the estimation of the speed of convergence of $\tau_{\mu,j}$
towards $\tau_\mu$.  Eventually, Section~\ref{ubi} is devoted to the
version of Theorem~\ref{cascades_ren2} needed to get (\ref{DIM}).

The techniques presented in this paper can be applied to derive
similar results for other statistically self-similar
$[0,1]$-martingales described in \cite{BM1,BACRYMUZY,BM2}.

%%%%%%%%%%%%%%%%%%%%%%%%%%%%%%%%%%%%%%%%%%%%%%%%
%%%%%%%%%%%%%%%%%%%%%%%%%%%%%%%%%%%%%%%%%%%%%%%%

%%%%%%%%%%%%%%%%%%%%%%%%%%%%%%%%%%%%%%%%%%%%%%%%
%%%%%%%%%%%%%%%%%%%%%%%%%%%%%%%%%%%%%%%%%%%%%%%%
%%%%%%%%%%%%%%%%%%%%%%%%%%%%%%%%%%%%%%%%%%%%%%%%
%%%%%%%%%%%%%%%%%%%%%%%%%%%%%%%%%%%%%%%%%%%%%%%%
\section{General estimates for the growth speed 
in singularity sets}
\label{principes}
%%%%%%%%%%%%%%%%%%%%%%%%%%%%%%%%%%%%%%%%%%%%%%%%
%%%%%%%%%%%%%%%%%%%%%%%%%%%%%%%%%%%%%%%%%%%%%%%%
\subsection{Measure of fine level sets: 
a neighboring boxes condition}\label{princ}

Let $(\Omega,\mathcal{B},\mathbb{P})$ stand for the probability space
on which the random variables in this paper are defined. Fix an
integer $b\ge 2$.

Let $\mathcal{A}=\{0,\dots,b-1\}$. For every $w\in
\mathcal{A}^*=\bigcup_{j\ge 0}\mathcal{A}^j$
($\mathcal{A}^0:=\{\emptyset \}$), let $I_w$ be the closed $b$-adic
subinterval of $[0,1]$ naturally encoded by $w$. If $w\in
\mathcal{A}^j$, we set $|w|=j$ .

For $n\ge 1$ and $0\le k\le b^n-1$, $I_{n,k}$ denotes the interval
$[kb^{-n},(k+1)b^{-n})$. If $t\in [0,1)$, $k_{n,t}$ is the unique
integer such that $t\in [k_{n,t}b^{-n},(k_{n,t}+1)b^{-n})$. We denote
by $w^{(n)}(t)$ the unique element $w$ of $\mathcal{A}^n$ such that
$I_w=[k_{n,t}b^{-n},(k_{n,t}+1)b^{-n}]$.

With $w\in \mathcal{A}^j$ can be associated a unique number $i(w)\in
\{0,1,\ldots, b^j-1\}$ such that
$I_w=[i(w)b^{-j},(i(w)+1)b^{-j}]$. Then, if $(v,w)\in \mathcal{A}^j$,
$\delta(v,w)$ stands for $|i(v)-i(w)|$. 

\smallskip

Let $\mu$ and $m$ be two positive Borel measures with supports equal
to $[0,1]$. 

Let $\widetilde\varepsilon=(\varepsilon_n)_{n\ge 0}$ be a positive
sequence, $N\ge 1$, and $\alpha\ge 0$.  

We consider a slight refinement of the sets introduced in (\ref{wtE}):
For $p\ge 1$,  define
\begin{eqnarray}
 \!\! \!\! E^{\mu}_{\alpha,p}(N,\widetilde\varepsilon) \!\!  & \!\!  =
\!\!  & \!\! \left \{ t\in [0,1]\!:\!
\begin{cases} \forall \, n\ge p, \ \forall w\in
\mathcal{A}^n \mbox{ such that }
\delta(w,w^{(n)}(t))\le N,\\ \ \forall\ \gamma\in\{-1,1\}, \mbox{ one
has } b^{\gamma n(\alpha-\gamma\varepsilon_n)}\mu(I_w)^{\gamma}\le 1
\end{cases} \!\!\! \!\!\!\right \}\\
\label{defebm}
 \!\!\!\!  \!\! \mbox{and } E^{\mu}_\alpha(N,\widetilde\varepsilon)
 \!\!  & \!\!  = \!\!  & \!\!  \bigcup_{p\ge
 1}E^{\mu}_{\alpha,p}(N,\widetilde\varepsilon).
\end{eqnarray}
This set contains the points $t$ for which, at each scale $n$ large
enough, the $\mu$-measure of the $2N+1$ neighbors of $I_{n,k_t}$
belongs to $[b^{-
n(\beta+\varepsilon_n)},b^{-j(\beta-\varepsilon_n)}]$. The information
on neighboring intervals is involved in proving (\ref{DIM}).

For $n\ge 1$ and $\varepsilon,\eta>0$, let us define the quantity
\begin{equation}
\label{defsneej}
S^{N,\varepsilon,\eta}_n(m,\mu,\alpha)= \sum_{\gamma\in
\{-1,1\}}b^{n(\alpha-\gamma\varepsilon)\gamma\eta} \hspace{-2mm}
\sum_{v,w\in \mathcal{A}^n: \ \delta(v,w)\le N} \hspace{-2mm}
m(I_v)\mu(I_w)^{\gamma\eta}.
\end{equation}
The following result is already established in \cite{BSGIBBS}, but we
give the proof for completeness.
%%%%%%%%%%%%%%%%%%%%%%%%%%%%%%%%%%%%%
\begin{proposition}\label{principe1}
Let $(\eta_n)_{n\ge 1}$ be a positive sequence. \\ If $\sum_{n\ge
1}S^{N,\varepsilon_n,\eta_n}_n(m,\mu,\alpha)<+\infty$, then
$E^{\mu}_\alpha(N,\widetilde\varepsilon )$ is of full
$m$-measure.
\end{proposition} 

%%%%%%%%%%%%%%%%%%%%%%%%%%%%%%%%%%%%%%%%
\begin{remark}
Similar conditions were used in \cite{BBP} to obtain a comparison
    between the multifractal formalisms of \cite{BRMICHPEY} and
    \cite{OLSEN}.
\end{remark}
%%%%%%%%%%%%%%%%%%%%%%%%%%%%%%%%%%%%%%%%
%%%%%%%%%%%%%%%%%%%%%%%%%%%%%%%%%%%%%%%%%%%%%%%%
\begin{proof}
For $\gamma\in\{-1,1\}$ and $n\ge 1$, let us define
\begin{equation}\label{ii}
E^{\mu}_\alpha(N,\varepsilon_n,\gamma)=\left\{t\in
[0,1]:\begin{cases}\forall\ w\in\mathcal{A}^n \mbox{ such that }
\delta(w,w^{(n)}(t))\le N,\\ \mbox{ one has } b^{\gamma
n(\alpha-\gamma\varepsilon_n)}\mu(I_w)^{\gamma}\le
1\end{cases}\hspace{-4mm}\right \}.
\end{equation}

\noindent
For $t\in [0,1]$, if there exists (a necessarily unique)
$w\in\mathcal{A}^n$ such that $i(w)-i(w^{(n)}(t))=k$, this word $w$ is
denoted $w^{(n)}_k(t)$. For $\gamma\in \{-1,1\}$, let $S_{n,\gamma}
=\sum_{-N\le k\le N}m_k$ with
$$m_k=m\left (\left\{t\in [0,1]: i(w)-i(w^{(n)}(t))=k\Rightarrow
b^{\gamma n(\alpha-\gamma\varepsilon_n)}\mu(I_w)^{\gamma}>1 \right
\}\right ).
$$ 
One clearly has
\begin{equation}\label{BC}
m\Big ((E^{\mu}_\alpha (N,\varepsilon_n,-1))^c\bigcup E^{\mu}_\alpha
(N,\widetilde\varepsilon_n,1))^c\Big )\le S_{n,-1}+S_{n,1}.
\end{equation}
Fix $\eta_n>0$ and $-N\le k\le N$. Let $Y(t)$ be random variable
defined to be equal to $b^{\gamma n(\alpha-\gamma\varepsilon_n)\eta_n}
\mu\big (I_{w^{(n)}_k(t)}\big )^{\gamma \eta_n}$ if $w^{(n)}_k(t)$
exists or 0 otherwise.  The Markov inequality applied to $Y(t)$ with
respect to $m$ yields $m_k \le \int Y(t) dm(t)$. Since $Y$ is constant
over each $b$-adic interval $I_v$ of generation $n$, we get
\begin{eqnarray*}
 m_k \le\sum_{v,w \in \mathcal{A}^n: \, i(w)-i(v)=k}
b^{n(\alpha-\gamma\varepsilon_n)\gamma \eta_n}m(I_v) \mu(I_w)^{\gamma
\eta_n}.
\end{eqnarray*}
Summing over $|k|\le N$ yields $S_{n,-1}+S_{n,1}\le
S^{N,\varepsilon_n,\eta_n}_n(m,\mu,\alpha)$.  The conclusion follows
from (\ref{BC}) and from the Borel-Cantelli Lemma.
\end{proof}
%%%%%%%%%%%%%%%%%%%%%%%%%%%%%%%%%%%%%

%%%%%%%%%%%%%%%%%%%%%%%%%%%%%%%%%%%%%
%%%%%%%%%%%%%%%%%%%%%%%%%%%%%%%%%%%%%
\subsection{Uniform growth speed in singularity sets}\label{princ2}

Let $\Lambda$ be a set of indexes, and $\Omega^*$ a measurable subset
of $\Omega$ of probability 1. Some notations and technical assumptions
are needed to state the general result that we shall apply in
Section~\ref{CASCADES}. These assumptions describe a common situation
in multifractal analysis. In particular the measures in the following
sections satisfy these requirements.

%%%%%%%%%%%%%%%%%%%%%%%%%%%%%%%%%%%%%
\smallskip

$\bullet$ For every $\omega\in\Omega^*$, we consider two sequences of
families of measures, namely $\Big (\{\mu^w_{\lambda}\}_{\lambda\in
\Lambda}\Big)_{w\in \mathcal{A}^*}$ and
$\Big(\{m^w_{\lambda}\}_{\lambda\in \Lambda}\Big)_{w\in
\mathcal{A}^*}$ (indexed by $\mathcal{A}^*$) such that for every $w\in
\mathcal{A}^*$, the elements of the families
$\{\mu^w_{\lambda}\}_{\lambda\in \Lambda}$ and
$\{m^w_{\lambda}\}_{\lambda\in \Lambda}$ are positive finite Borel
measures whose support is $[0,1]$. For $\nu\in\{\mu,m\}$,
$\{\nu^\emptyset _{\lambda}\}_{\lambda\in \Lambda}$ is written
$\{\nu_{\lambda}\}_{\lambda\in \Lambda}$.

%%%%%%%%%%%%%%%%%%%%%%%%%%%%%%%%%%%%%
\smallskip

$\bullet$ We consider an integer $N\ge 1$, a positive sequence
$\widetilde\varepsilon=(\varepsilon_{n})_{n\ge 1}$, and a family of
positive numbers $(\alpha_\lambda)_{\lambda\in\Lambda}$. Then,
remembering (\ref{ii}) let us consider for every $j\ge 0$,
$w\in\mathcal{A}^j$ and $p\ge 1$ the sets
\begin{equation}
\label{deflevel}
E^{\mu^{w}_\lambda}_{\alpha_\lambda,p}(
N,\widetilde\varepsilon)=\bigcap_{n\ge p
}E_{\alpha_\lambda}^{\mu^{w}_\lambda}(N, \varepsilon_{n},-1)\cap
E_{\alpha_\lambda}^{\mu^{w}_\lambda}(N, \varepsilon_{n},1).
\end{equation}
The sets $\{ E^{\mu^w_\lambda}_{\alpha_\lambda,n}(
N,\widetilde\varepsilon)\}_n$ form a non-decreasing sequence. We
assume that $m_\lambda^{w}$ is concentrated on $\lim_{p\rightarrow
+\infty}
E^{\mu^{w}_\lambda}_{\alpha_\lambda,p}(N,\widetilde\varepsilon)$. One
defines the growth speed of $ E^{\mu^{w}_\lambda}_{\alpha_\lambda,p}(
N,\widetilde\varepsilon)$~as
\begin{equation}
\label{defn}
GS(m^{w}_\lambda,\mu^{w}_\lambda,\alpha_\lambda,
N,\widetilde\varepsilon)= \inf \left\{ p\ge 1 : m_\lambda^{w}\big
(E^{\mu^{w}_\lambda}_{\alpha_\lambda,p}( N,\widetilde\varepsilon)\big )
\ge {1}/{2} \right\}.
\end{equation}
This number, which may be infinite, is a measure of the number $p$ of
generations needed for $E^{\mu^{w}_\lambda}_{\alpha_\lambda,p}(
N,\widetilde\varepsilon)$ to recover a certain given fraction (here
chosen equal to 1/2) of the measure $m_\lambda^{w}$. Since
$\mu_\lambda^w(\lim_ {p\rightarrow +\infty}
E^{\mu^{w}_\lambda}_{\alpha_\lambda,p }(N,\widetilde\varepsilon))=1$,
$GS(m^{w}_\lambda,\mu^{w}_\lambda,\alpha_\lambda,N,\widetilde\varepsilon)
<\infty$.

%%%%%%%%%%%%%%%%%%%%%%%%%%%%%%%%%%%%%

%%%%%%%%%%%%%%%%%%%%%%%%%%%%%%%%%%%%%
\smallskip

$\bullet$ We assume that for every positive sequence $\widetilde\eta=
(\eta_j)_{j\ge 0}$, there exist a random vector $\big
(U(\widetilde\eta),V (\widetilde\eta)\big ) \in \mathbb{R}_+\times
\mathbb{R}_+^{\mathbb{N}}$ and a sequence
$(U^w,V^w)_{w\in\mathcal{A}^*}$ of copies of $\big
(U(\widetilde\eta),V (\widetilde\eta)\big )$ and finally a sequence
$(\psi_j(\widetilde\eta))_{j\ge 0}$, such that for $\mathbb{P}$-almost
every $\omega\in\Omega^*$, 
\begin{equation}
\label{bound}
\forall \, w\in \mathcal{A}^*, \, \forall n\ge \psi_{|w|}
(\widetilde\eta) ,  \begin{cases}
U^w\le \inf_{\lambda\in \Lambda}\Vert m^w_\lambda\Vert,\\
V^w_n \ge \sup_{\lambda\in\Lambda} S_n ^{N,\varepsilon_n,\eta_n}
(m^w_\lambda,\mu^w_\lambda,\alpha_\lambda),
\end{cases} \!\!
\end{equation}
where $S_n^{N,\varepsilon_{n},
\eta_n}(m^w_\lambda,\mu^w_\lambda,\alpha_\lambda)$ is defined in
(\ref{defsneej}) (if $w\in {\mathcal A}_j$, remember that
$|w|=j$). This provides us with a uniform control over
$\lambda \in \Lambda$ of the families of measures
$(m^w_\lambda,\mu^w_\lambda)_{w\in \mathcal{A}^*}$.

%%%%%%%%%%%%%%%%%%%%%%%%%%%%%%%%%%%%%

%%%%%%%%%%%%%%%%%%%%%%%%%%%%%%%%%%%%%
%%%%%%%%%%%%%%%%%%%%%%%%%%%%%%%%%%%%%
\begin{proposition}[{\bf Uniform growth speed in singularity sets}]\label{ren}
Assume that two sequences of positive numbers $\widetilde
\eta=(\eta_j)_{j\ge0} $ and $(\rho_j)_{j\ge 0}$ are fixed.\\ Let
$(\mathcal{S}_j)_{j\ge 0}$ be a sequence of integers such that
$\mathcal{S}_j\ge \psi_j(\widetilde\eta)$. If
\begin{equation}
\label{C1}
\sum_{j\ge 0}b^j\mathbb{P}\Big (U(\widetilde\eta) \le b^{-\rho_j}\Big
)<\infty \mbox{ and }\sum_{j\ge 0}b^jb^{\rho_j} \sum_{k\ge
\mathcal{S}_j} \mathbb{E} \Big (V_k(\widetilde\eta)\Big )<\infty,\!\!
\end{equation}
then, with probability one, for $j$ large enough, for every $w\in
\mathcal{A}^j$ and $\lambda\in \Lambda$, one has
$GS(m^{w}_\lambda,\mu^{w}_\lambda,\alpha_\lambda,N,\widetilde\varepsilon)\le
\mathcal{S}_j$.
\end{proposition}
%%%%%%%%%%%%%%%%%%%%%%%%%%%%%%%%%%%%%
%%%%%%%%%%%%%%%%%%%%%%%%%%%%%%%%%%%%%

%%%%%%%%%%%%%%%%%%%%%%%%%%%%%%%%%%%%%
%%%%%%%%%%%%%%%%%%%%%%%%%%%%%%%%%%%%%
\begin{proof}
Fix $j\ge 1$ and $w\in \mathcal{A}^j$. As shown in the proof of
Proposition \ref{principe1}, for every $n\ge 1$ and every
$\lambda\in\Lambda$, one can write
\begin{eqnarray*}
m^w_\lambda\Big (\big(E^{\mu^w_\lambda}_{\alpha_\lambda}(N,
\varepsilon_{n},-1)\big )^c\cup \big
(E^{\mu^w_\lambda}_{\alpha_\lambda}(N, \varepsilon_{n},1)\big )^c\Big
) \le S_n^{N,\varepsilon_{n},\eta_n}(
m^w_\lambda,\mu^w_\lambda,\alpha_\lambda).
\end{eqnarray*}
Thus, using (\ref{bound}), one gets
\begin{equation}
\label{ineg1}
m^w_\lambda \Big ( \bigcup_{n\ge \mathcal{S}_j} \big
(E^{\mu^w_\lambda}_{\alpha_\lambda}(N, \varepsilon_{n},-1)\big )^c\cup
\big (E^{\mu^w_\lambda}_{\alpha_\lambda}(N, \varepsilon_{n},1)\big
)^c\Big )\le \sum_{n\ge \mathcal{S}_j}V^w_n.\!\!
\end{equation}
 Let us apply the ``statistical self-similar control'' (\ref{C1})
combined with the Borel-Cantelli lemma.  On the one hand, the left
part of (\ref{C1}) yields $ \sum_{j\ge 1}\mathbb{P}\Big (\exists \
w\in \mathcal{A}^j,\ U^w\le b^{-\rho_j}\Big )<\infty$. Hence, with
probability one, for $j$ large enough and for all $w\in\mathcal{A}^j$,
\begin{equation}
\label{aa}
\sup_{\lambda\in \Lambda}\Vert m^w_\lambda\Vert^{-1}\le (U^w)^{-1}\le
b^{\rho_j}.
\end{equation}
On the other hand, the right part of (\ref{C1}) yields
\begin{eqnarray*}
\sum_{j\ge 1}\mathbb{P}\Big (\exists \ w\in \mathcal{A}^j,\ b^{ \rho_j
}\sum_{n\ge \mathcal{S}_j}V^w_n\ge 1/2\Big ) \le 2\sum_{j\ge
1}b^jb^{\rho_j }\mathbb{E}\Big (\sum_{n\ge \mathcal{S}_j}V^w_n \Big
)<\infty.
\end{eqnarray*}
This implies that with probability one, $ b^{\rho_j}\sum_{n\ge
\mathcal{S}_j}V^w_n<1/2 $ for every $j$ large enough and for all
$w\in\mathcal{A}^j$.  

Thus, by (\ref{aa}), $\sup_{\lambda\in \Lambda}\frac{\sum_{n\ge
\mathcal{S}_j}V^w_n}{\Vert m_\lambda^w\Vert }<1/2$.  Combining this
with (\ref{ineg1}) and (\ref{defn}), one gets that for every
$\lambda\in \Lambda$, $GS(m^w_\lambda,\mu^w_\lambda,\alpha_\lambda,
N,\widetilde\varepsilon)\le \mathcal{S}_j$.
\end{proof}
%%%%%%%%%%%%%%%%%%%%%%%%%%%%%%%%%%%%%%

%%%%%%%%%%%%%%%%%%%%%%%%%%%%%%%%%%%%%%
%%%%%%%%%%%%%%%%%%%%%%%%%%%%%%%%%%%%%%
%%%%%%%%%%%%%%%%%%%%%%%%%%%%%%%%%%%%%%
%%%%%%%%%%%%%%%%%%%%%%%%%%%%%%%%%%%%%%

%%%%%%%%%%%%%%%%%%%%%%%%%%%%%%%%%%%%%
\section{Main results for independent random cascades} 
\label{CASCADES}
%%%%%%%%%%%%%%%%%%%%%%%%%%%%%%%%%%%%%
%%%%%%%%%%%%%%%%%%%%%%%%%%%%%%%%%%%%%
\subsection{Definition}\label{IRC}

Let $v=(v_1,\ldots, v_{|v|})\in\mathcal{A}^*$. We need the following
truncation notation: For every $k\in \{1,\ldots,|v|\} $, $v|k$ is the
word $(v_1,\ldots, v_{k}) \in \mathcal{A}^{k}$, and by convention
$v|0$ is the empty word~$\emptyset$.

We focus in this paper on the measures introduced in \cite{M} and more
recently in \cite{B2}. A measure $\mu(\omega)$ is said to be an
independent random cascade if it has the following property: There
exist a sequence of random positive vectors $\big (W(w)=(W_0(w),\dots,
W_{b-1}(w))\big )_{w\in \mathcal{A}^{*}}$ and a sequence of random
measures $(\mu^{w})_{w\in \mathcal{A}^*}$ such that

\smallskip
\noindent
$\bullet$ {\bf (P1)} for all $v,w\in\mathcal{A}^*$, $\mu (I_{vw})=
\mu^{v}(I_w)\prod_{k=0}^{|v|-1}W_{v_{k+1}}(v|k)$ ($\mu^{\emptyset} =
\mu$),

\smallskip
\noindent
$\bullet$
{\bf (P2)} the random vectors $W(w)$, for $w\in \mathcal{A}^*$, are
i.i.d. with a vector $W=(W_0,\ldots,W_{b-1})$ such that
$\sum_{k=0}^{b-1}\mathbb{E}(W_k)=1$,

\smallskip
\noindent
$\bullet$
{\bf (P3)} for all $v\in \mathcal{A}^* $, $\big (\mu^{v}(I_w)\big
)_{w\in \mathcal{A}^*}\equiv \big (\mu(I_w)\big )_{w\in
\mathcal{A}^*}$. Moreover, for every $j\ge 1$, the measures
$\mu^{v}$, $v\in \mathcal{A}^j $, are mutually independent,

\smallskip
\noindent
$\bullet$
{\bf (P4)} for every $j\ge 1$, the $\sigma$-algebras $\sigma\big
(W(w):\ w\in \cup_{0\le k\le j-1}\mathcal{A}^{k}\big )$ and $\sigma
(\mu^v(I_w):\ v\in \mathcal{A}^{j},\ w\in\mathcal{A}^*)$ are
independent.

\smallskip

Let $\big (W(w)\big )_{w\in \mathcal{A}^*}$ be as above.  For
$q\in\mathbb{R}$ define the function 
\begin{equation}
\label{zero}
\widetilde\tau_\mu(q)=-\log_b\mathbb{E}\Big(\,\sum_{k=0}^{b-1}W_k^q\,
\Big)  \in\mathbb{R}\cup\{-\infty\}.
\end{equation}
The two classes of measures we deal with are the following.

\smallskip

%%%%%%%%%%%%%%%%%%%%%%%%%%%%%%%%%%%%%
%%%%%%%%%%%%%%%%%%%%%%%%%%%%%%%%%%%%%
\noindent{\bf Non-degenerate multiplicative martingales when
$\widetilde\tau'_\mu(1^-)>0$.} With probability one, $\forall \,
v\in\mathcal{A}^*$, the sequence of measures $\{\mu^v_j\}_{j\ge 0}$
\begin{equation}\label{def20}
\frac{d\mu^v_j}{d\ell}(t)=b^j\prod_{k=0}^{j-1}W_{t_{k+1}}(v\cdot t|k)
\end{equation} defined on $[0,1]$ converges weakly, as $n\to\infty$ 
to a measure $\mu^v$.  For $\mu=\mu^\emptyset$:\\ (1) Properties {\bf
(P1)} to {\bf (P4)} are satisfied;\\ (2) If
$\widetilde\tau'_\mu(1^-)>0$, the total masses $\Vert\mu^v\Vert$ are
almost surely positive, and their expectation is equal to 1
(\cite{KP,DL}).

\medskip

%%%%%%%%%%%%%%%%%%%%%%%%%%%%%%%%%%%%%
%%%%%%%%%%%%%%%%%%%%%%%%%%%%%%%%%%%%%
\noindent{\bf The modified construction in the critical case
$\widetilde \tau'_\mu(1^-)=0$}. Suppose that $\widetilde
\tau'_\mu(1^-)=0$ and $\widetilde\tau_\mu(h)>-\infty$ for some $h>1$.
Then, with probability one, for all $v\in\mathcal{A}^*$, the function
of $b$-adic intervals 
\begin{equation}
\label{def2}
\mu^v(I_w)=\lim_{j\to\infty}-\sum_{u\in\mathcal{A}^j}\Big
(\prod_{k=0}^{|w|+j-1}W_{(w\cdot u)_{k+1}}\big (v\cdot (w\cdot
u|k)\big )\Big )\log \prod_{k=0}^{|w|+j-1}W_{(w\cdot u)_{k+1}}\big
(v\cdot (w\cdot u|k)\big )
\end{equation}
is well defined and yields a positive Borel measure whose support is
$[0,1]$ (\cite{B2,L2}). For $\mu=\mu^\emptyset$,\\ (1) properties {\bf
(P1)} to {\bf (P4)} are satisfied;\\ (2)
$\mathbb{E}(\Vert\mu\Vert^h)<\infty$ for $h\in [0,1)$ but
$\mathbb{E}(\Vert\mu\Vert)=\infty$.

%%%%%%%%%%%%%%%%%%%%%%%%%%%%%%%%%%%%%
%%%%%%%%%%%%%%%%%%%%%%%%%%%%%%%%%%%%%
\subsection{Analyzing measures} 
In both above cases we define $J$ as the interior of the interval
$\{q\in\mathbb{R}: \widetilde
\tau'_\mu(q)q-\widetilde\tau_\mu(q)>0\}$. One always has $(0,1)\subset
J$ and $J\subset (-\infty,1)$ if $\widetilde\tau'_\mu(1^-)=0$. We
assume that:\\ - If $\widetilde\tau'_\mu(1^-)>0$ , $J$ contains the
closed interval $[0,1]$,\\ - If $\widetilde\tau'_\mu(1^-)=0$ then
$0\in J$.

\smallskip

For $q\in J$, $v\in\mathcal{A}^*$, $j\ge 1$, let $\mu^v_{q,j}$ be the
measure defined as $\mu^v_j$ in (\ref{def20}) but with the sequence
$\big (W_q(v\cdot w)=(b^{\widetilde\tau_\mu(q)}W_0(v\cdot
w)^q,\dots,b^{\widetilde \tau_\mu(q)}W_{b-1}(v\cdot w)^q)\big )_w$
instead of $\big (W(v\cdot w)\big )_w$. It is proved in \cite{B2} that
there exists a subset $\Omega^*$ of $\Omega$ of probability 1 such
that  $\forall \omega\in\Omega^*$, $\forall \,v\in\mathcal{A}^*$ and
$\forall \,q\in J$, the sequence $\mu^v_{q,j}$ converges weakly to a
positive measure $\mu^v_q$.

\smallskip

If one denotes $\mu_q^\emptyset=\mu_q$, $Y_q =
\Vert\mu^\emptyset_q\Vert$, and $Y_q(v) = \Vert\mu^v_q\Vert$ for
$v\in\mathcal{A}^*$, it is proved in \cite{Biggins,B2} that with
probability one the mappings $q\in J\mapsto Y_q(v)$ are analytic and
positive. Moreover, \cite{B2} proves that
$\tau_\mu\equiv\widetilde\tau_\mu$ on $J$ almost surely.
 
Eventually, one can see that $J\supset\mathbb{R}_+$
(resp. $\mathbb{R}_-$) if and only if
$\widetilde\tau_\mu(hq)-h\widetilde\tau_\mu(q)>0$ for all $q\in
\mathbb{R}_+$ (resp. $\mathbb{R}_-$) and $h>1$, which amounts to
saying that $\forall \,q\in \mathbb{R}_+$ (resp. $\mathbb{R}_-$),
$\mathbb{E}(Y_q^h)<\infty$ and $h>1$ (see the proof of
Lemma~\ref{Momentspos}).

%%%%%%%%%%%%%%%%%%%%%%%%%%%%%%%%%%%%%
%%%%%%%%%%%%%%%%%%%%%%%%%%%%%%%%%%%%%
\subsection{Main results}

 Recall that for an independent random cascade $\mu$, if
$\widetilde\tau'_\mu(1^-)>0$ we assume that $J$ contains $[0,1]$, and
if $\widetilde\tau'_\mu(1^-)=0$ then $J\subset (-\infty,1)$, and we
assume that $0\in J$.
%%%%%%%%%%%%%%%%%%%%%%%%%%%%%%%%%%%%%
%%%%%%%%%%%%%%%%%%%%%%%%%%%%%%%%%%%%%
\begin{theorem}\label{cascades}
Let $\mu$ be an independent random cascade. Let $N$ be an integer $\ge
 1$ and $\widetilde\varepsilon =(\varepsilon_n)_{n\ge 1}$ a sequence
 of positive numbers going to 0. Assume that $\forall \alpha>0$, the
 series $\sum_{n\ge 1} n{b}^{-n\alpha\varepsilon^2_n }$ converges.

Then, with probability one, for every $q\in J$,
 $\tau_\mu(q)=\widetilde\tau_\mu(q)$, and the two level sets
 $E^{\mu_q}_{\tau_\mu'(q)q-\tau_\mu(q)}(N,\widetilde \varepsilon)$ and
 $E^{\mu}_{\tau_\mu'(q)}(N,\widetilde\varepsilon)$ are both of full
 $\mu_q$-measure.
\end{theorem}
%%%%%%%%%%%%%%%%%%%%%%%%%%%%%%%%%%%%%%%%
%%%%%%%%%%%%%%%%%%%%%%%%%%%%%%%%%%%%%%%%

%%%%%%%%%%%%%%%%%%%%%%%%%%%%%%%%%%%%%
%%%%%%%%%%%%%%%%%%%%%%%%%%%%%%%%%%%%%
\begin{remark} 
The conclusions of Theorem~\ref{cascades} hold as soon as
\begin{equation}
\label{assump}
\exists \eta>0 \mbox{ such that for every $n$, }\varepsilon_n\ge
n^{-1/2}\log(n)^{{1}/{2}+\eta}.
\end{equation}
\end{remark}
%%%%%%%%%%%%%%%%%%%%%%%%%%%%%%%%%%%%%
%%%%%%%%%%%%%%%%%%%%%%%%%%%%%%%%%%%%%

%%%%%%%%%%%%%%%%%%%%%%%%%%%%%%%%%%%%%%%%
%%%%%%%%%%%%%%%%%%%%%%%%%%%%%%%%%%%%%%%%
\begin{theorem}[{\bf Growth speed in H\"older singularity sets}]
  \label{cascades_ren} 
Under the assumptions of Theorem \ref{cascades}, assume that
 $(\varepsilon_n)_n$ satisfies (\ref{assump}) and that there exists
 $A>1$ such that with probability one $A^{-1}\le W_i$ (resp. $W_i\le
 A$) for all $0\le i\le b-1$ . Let $K$ be a compact subinterval of
 $J\cap \mathbb{R}_+$ (resp. $J\cap \mathbb{R}_-$).

Then, with probability one, for $j$ large enough, for all $q\in K$ and
$w\in \mathcal{A}^j$, 
$$\max \left (GS\big(\mu_q^w,\mu^w,\tau_\mu'(q), N,\widetilde
\varepsilon\big ),GS\big(\mu_q^w,\mu_q^w,\tau_\mu'(q)q-\tau_\mu(q),
N,\widetilde \varepsilon\big) \right)\le \mathcal{S}_j$$ with
$\mathcal{S}_j=\big [ \exp \big (\big (j\log
(j)^{\eta}\big)^{\frac{1}{1+2\eta}}\big )\big ] $.

If there exists $\eta>0$ such that for every $n$, $\varepsilon_n\geq
\log (n)^{-\eta}$, the above conclusion holds with $\mathcal{S}_j=\big
[ j\log (j)^{\eta'}\big ]$, for any $\eta'>2\eta$.

\end{theorem}
%%%%%%%%%%%%%%%%%%%%%%%%%%%%%%%%%%%%%
%%%%%%%%%%%%%%%%%%%%%%%%%%%%%%%%%%%%%%%%
\begin{remark}
In Theorem~\ref{cascades_ren}, the first choice of $\varepsilon_n$
corresponds to the ``best'' choice for the speed of convergence
$\varepsilon_n$, and the growth speed $\mathcal{S}_j$ is very slow. To
the contrary, the second choice for $\varepsilon_n$ is ``worse'', but
as a counterpart $\mathcal{S}_j$ is improved.

We assume that the number of neighbors $N$ is fixed. In fact, it is
not difficult to consider a sequence of neighbors $N_n$ simultaneously
with the speed of convergence $\varepsilon_n$. This number $N_n$ can
then go to $\infty$ under the condition that $\log
N_n=o(n\varepsilon_n^2)$.  Another modification would consist in
replacing the fixed fraction $f$ in (\ref{GS}) by a fraction $f_j$
going to 1 as $j$ goes to $\infty$. The choice $f_j=1-b^{-s_j}$ with
$s_j=j$ is convenient. These two improvements yield technical
complications, but comparable results are easily derived from the
proofs we propose.
\end{remark}
%%%%%%%%%%%%%%%%%%%%%%%%%%%%%%%%%%%%%

The growth speed obtained in Theorem~\ref{cascades_ren} can be
improved by considering results valid only almost surely, for almost
every $q$, $\mu_q$ almost-everywhere.

Recall that if $t\in [0,1)$ and $j\ge 1$, $w^{(j)}(t)$ is the unique
element $w$ of $\mathcal{A}^j$ such that $t\in
[i(w)b^{-j},(i(w)+1)b^{-j})$.

%%%%%%%%%%%%%%%%%%%%%%%%%%%%%%%%%%%%%
%%%%%%%%%%%%%%%%%%%%%%%%%%%%%%%%%%%%%
\begin{theorem}[{\bf Improved growth speed}]
\label{cascades_ren2}
Under the assumptions of Theorem \ref{cascades}, fix $\kappa>0$ and
assume that (\ref{assump}) holds. For $j\ge 2$, let $\mathcal{S}_j=\big
[j\log (j)^{-\kappa} \big ]$.

\smallskip
\noindent (1) For every $q\in J$, with probability one, the property
$\mathcal{P}(q)$ holds, where $\mathcal{P}(q)$ is: For $\mu_q$-almost
every $t\in [0,1)$, if $j$ is large enough, for $w=w^{(j)}(t)$ one has
$$
\max \left (GS\big(\mu_q^w,\mu^w,\tau_\mu'(q), N,\widetilde
\varepsilon\big ),GS\big(\mu_q^w,\mu_q^w,\tau_\mu'(q)q-\tau_\mu(q),
N,\widetilde \varepsilon\big) \right)\le \mathcal{S}_j.$$

\noindent (2) With probability one, for almost every $q\in J$,
$\mathcal{P}(q)$ holds.
\end{theorem}
%%%%%%%%%%%%%%%%%%%%%%%%%%%%%%%%%%%%%
%%%%%%%%%%%%%%%%%%%%%%%%%%%%%%%%%%%%%

For $w\in\mathcal{A}^*$, $n\ge 1$ and $q\in J$, let 
\begin{equation}
\label{defnld}
\mathcal{N}_n(\mu^w,\alpha,\varepsilon_n)\!= \!\#\big\{b\mbox{-adic
  box }I \mbox{ of scale }n\!:\!|I|^{\alpha+\varepsilon_n}\!\le \mu^w
  (I )\!\le |I|^{\alpha-\varepsilon_n}\big\}.
\end{equation}
Remember that $\tau_\mu=\widetilde\tau_\mu$ on~$J$.

%%%%%%%%%%%%%%%%%%%%%%%%%%%%%%%%%%%%%%%%
%%%%%%%%%%%%%%%%%%%%%%%%%%%%%%%%%%%%%
\begin{theorem}[{\bf Renewal speed of large
deviations spectrum}]\label{cascades_dev} Under the assumptions of
 Theorem \ref{cascades}, let us also assume that (\ref{assump}) holds,
 $J=\mathbb{R}$ (in particular $\widetilde \tau_\mu'(1)>0$) and there
 exists $A>1$ such that with probability one $A^{-1}\le W_i\le A$ for
 all $0\le i\le b-1$. Let $\mathcal{S}_j$ be defined as in
 Theorem~\ref{cascades_ren}. Let $K$ be a compact subinterval of
 $\mathbb{R}$ and $\beta=1+\max_{ q \in K}|q|$.

Then, with probability one, for $j$ large enough, for all $q\in K$ and
$w\in \mathcal{A}^j$, and for all $n\ge \mathcal{S}_j$, one has
$$
Y_q^wb^{n(\widetilde
\tau_\mu'(q)q-\widetilde\tau_\mu(q)-\beta\varepsilon_n)}\le \mathcal
N_n\big (\mu^{w},\widetilde \tau_\mu'(q),\varepsilon_n\big )\le
Y_q^wb^{n(\widetilde
\tau_\mu'(q)q-\widetilde\tau_\mu(q)+\beta\varepsilon_n)}.
$$
\end{theorem}
%%%%%%%%%%%%%%%%%%%%%%%%%%%%%%%%%%%%%%%%

For $w\in \mathcal{A}^*$, $n\ge 1$ and $q\in\mathbb{R}$, let us
introduce the functions ($\tau^\emptyset_n$ associated with
$\mu^\emptyset=\mu$ is simply denoted by $\tau_n$)
$$
\tau^w_n(q)=-\frac{1}{n}\log_b\sum_{v\in \mathcal{A}^n}\mu^w(I_v)^q.
$$

%%%%%%%%%%%%%%%%%%%%%%%%%%%%%%%%%%%%%

The speed of convergence obtained in  Theorem~\ref{convergence}
provides precisions on the estimator of the function $\tau_\mu$
discussed in \cite{colletkoukiou,OW}.

%%%%%%%%%%%%%%%%%%%%%%%%%%%%%%%%%%%%%%%%%%
\begin{theorem}[{\bf Convergence speed of $\tau^w_n$ toward
$\widetilde\tau_\mu$}]\label{cvtau}
\label{convergence} Under the assumptions of
Theorem~\ref{cascades_dev}, let $K$ be a compact subinterval of
$\mathbb{R}$. There exists $\theta>0$ and $\delta \in (0,1)$ such
that, with probability one,\\
1. for $j$ large enough, 
$\displaystyle
\big |\widetilde \tau_\mu(q)-\tau_j(q)\big |\le {\big |\log_b Y_q\big
|}\, {j^{-1}}+\theta \,{\log (j)}j^{-1};
$

\noindent
2.  for $j$ large enough, for every $n\ge j^{\delta}$, for every $
w\in \mathcal{A}^j$, $\displaystyle \big
|\widetilde\tau_\mu(q)-\tau^w_n(q)\big |\le {\big |\log_b Y^w_q\big
  |}\, {n^{-1}}+\theta \, {\log (n)}{n^{-1}}, $ with $|\log Y^w_q\big|
\le \theta \log (j)$.
\end{theorem}
%%%%%%%%%%%%%%%%%%%%%%%%%%%%%%%%%%%%%%%%%%
%%%%%%%%%%%%%%%%%%%%%%%%%%%%%%%%%%%%%%%%%%%

%%%%%%%%%%%%%%%%%%%%%%%%%%%%%%%%%%%%%%%%%%%
%%%%%%%%%%%%%%%%%%%%%%%%%%%%%%%%%%%%%%%%%%%
\subsection{Proof of Theorem~\ref{cascades}} Fix $K$ a compact subinterval
  of $J$ and $ \widetilde \eta = (\eta_n)_{n\ge 1}$ a bounded positive
sequence to be specified later. For $\omega\in\Omega^*$ and $q\in K$,
let us introduce (recall (\ref{defsneej}))
\begin{equation}
\label{deffj}
F_n(q)=S_n^{N,\varepsilon_n,\eta_n}\big (\mu_q,\mu,\tau_\mu'(q)\big )
\mbox{ and } G_n(q)=S_n^{N,\varepsilon_n,\eta_n}\big
(\mu_q,\mu_q,\tau_\mu'(q)q-\tau_\mu(q)\big ).
\end{equation}
We begin by giving estimates for $\mathbb{E}\big (H_n(q)\big )$ and
$\mathbb{E}\big (H_n'(q)\big )$ for $H\in \{F,G\}$.
%%%%%%%%%%%%%%%%%%%%%%%%%%%%%%%%%%%%%
\begin{lemma}\label{bb}
Under the assumptions of Theorem~\ref{cascades}, if
$\Vert\widetilde\eta\Vert_\infty$ is small enough, there exists
$C_K>0$ such that if $H\in \{F,G\}$,
\begin{equation}
\label{eq1}
 \forall \, q\in K, \ \max \big(\mathbb{E}\big (H_n(q)\big
),\mathbb{E}\big (H'_n(q)\big )\big )\le C_Knb^{-n 
(\varepsilon_n\eta_n+O(\eta_n^2)  )},
\end{equation}
where $O(\eta_n^2)$ is uniform over $q\in K$.
\end{lemma}
%%%%%%%%%%%%%%%%%%%%%%%%%%%%%%%%%%%%%
The proof of this lemma is postponed to the next subsection. 

Let $q_0$ be the left end point of $K$. Since $ \sup_{q\in
K}H_n(q)\le H_n(q_0)+\int_K|H_n'(q)|\, dq$, one has
\begin{equation}\label{crucial}
\mathbb{E}\big (\sup_{q\in K}H_n(q)\big )\le C_K(1+|K|) nb^{-n
(\varepsilon_n\eta_n+O(\eta_n^2) )}.
\end{equation}
Choosing $\eta_n={\varepsilon_n}/{A}$ with $A$ large enough yields
$\varepsilon_n\eta_n+O(\eta_n^2)\ge {\varepsilon_n^2}/{2A}$. Using 
the assumptions of Theorem~\ref{cascades},  we get the almost sure
convergence of $\sum_{n\ge 1}\sup_{q\in K}H_n(q)$ for $H\in
\{F,G\}$. We conclude with Proposition~\ref{principe1}.

%%%%%%%%%%%%%%%%%%%%%%%%%%%%%%%%%%%%%%%%%%%
%%%%%%%%%%%%%%%%%%%%%%%%%%%%%%%%%%%%%%%%%%%
\subsection{Proof of Lemma~\ref{bb}}
{\bf $\bullet$ The case $H=F$:} For $v,w\in\mathcal{A}^n$, $q\in J$
and $\gamma\in \{-1,1\}$, one can write
$$ 
\mu_q(I_v)\mu
(I_w)^{\gamma\eta_n}=Y_q(v)Y_1(w)^{\gamma\eta_n}b^{n\tau_\mu(q)}
\prod_{k=0}^{n-1}W_{v_{k+1}}(v|k)^qW_{w_{k+1}}(w|k)^{\gamma\eta_n}.
$$
Moreover, it follows from estimates of \cite{B2} that for
$\Vert\widetilde \eta\Vert_\infty$ small enough, the quantities
\begin{eqnarray*}
C'_K(\widetilde\eta) & = & \sup_{\substack{ q\in K, \ \gamma\in\{-1,1\}\\\
n\geq 1, \ v,w\in {\mathcal A}^n}}\Big ({\mathbb E}\Big
(\Big|\frac{d}{dq}Y_q(v) Y_1(w)^{\gamma\eta_n}\Big|\Big ) +{\mathbb
E}\big (Y_q(v)Y_1(w)^{\gamma\eta_n}\big )\Big )\\ \mbox{and }\ \
C''_K(\widetilde\eta) & = & \sup_{\substack{ q\in K,\ \gamma\in
\{-1,1\}\\\ n\geq 1,\ v,w\in {\mathcal A}^n,\ 0\le k\le
n-1}}\frac{{\mathbb E}\big
(\big|\frac{d}{dq}W_{v_{k+1}}(v|k)^qW_{w_{k+1}}(w|k)^{\gamma\eta_n}\big|\big
)}{{\mathbb E}\left
(W_{v_{k+1}}(v)^qW_{w_{k+1}}(w)^{\gamma\eta_n}\right )}
\end{eqnarray*}
are finite. Hence, due to the definition of $F_n(q)$ and the fact that
 $\widetilde\tau$ is continuously differentiable on $J$, there exists
 a constant $C_K(\widetilde\eta)$ such that for every $q\in K$, $\max
 \big(\mathbb{E}\big (F_n(q)\big ),\mathbb{E}\big (F'_n(q)\big )\big
 )\le C_K(\widetilde \eta)T_n(q), $ where
$$ T_n(q)=n \, b_n(q)\sum_{\substack{\gamma\in \{-1,1\},\\ v,w\in
\mathcal{A}^n,\ \delta(v,w)\le N}} \prod_{k=0}^{n-1}\mathbb{E}\big
(W_{v_{k+1}}(v|k)^qW_{w_{k+1}}(w|k)^{\gamma\eta_n}\big ),
$$ where $b_n(q)= b^{
n(\tau_\mu(q)+\gamma\eta_n(\tau_\mu'(q)-\gamma\varepsilon_n))}$. Let
us make the following important remark.

%%%%%%%%%%%%%%%%%%%%%%%%%%%%%%%%%%%%%%%%%
\begin{remark}\label{reduction} If
$v$ and~$w$ are words of length~$n$, and if $\bar{v}$ and $\bar{w}$
stand for their prefixes of length~$n-1$, then
$\delta(\bar{v},\bar{w})> k$ implies $\delta(v,w) >bk$. It implies
that, given two integers $n \ge m>0$ and two words $v$ and $w$ in
${\mathcal A}^n$ such that $b^{m-1}<\delta(v,w)\le b^m$, there are two
prefixes $\bar{v}$ and $\bar{w}$ of respectively $v$ and $w$ of common
length\ $n-m$\ such that $\delta(\bar{v},\bar{w})\le 1$. Moreover, for
these words $\bar{v}$ and $\bar{w}$, there are at most $b^{2m}$ pairs
$(v,w)$ of words in $\mathcal{A}^n$ such that $\bar{v}$ and $\bar{w}$
are respectively the prefixes of $v$ and $w$.
\end{remark} 
%%%%%%%%%%%%%%%%%%%%%%%%%%%%%%%%%%%%%%%%%

Due to Remark \ref{reduction} and the form of $ T_n(q)$, there exists
a constant $C'_K$ such that for all $q\in K$ and $n\ge 1$
$$ T_n(q)\le C'_K n b_n(q)\sum_{\substack{\gamma\in \{-1,1\},\\ v,w\in
\mathcal{A}^n,\ \delta(v,w)\le 1}} \prod_{k=0}^{n-1}\mathbb{E}\big
(W_{v_{k+1}}(v|k)^qW_{w_{k+1}}(w|k)^{\gamma\eta_n}\big ).
$$
 The situation is thus reducible to the case $N=1$. Now $T_n(q)\le
n\big (T_{n,1}(q)+T_{n,2}(q)\big ),$
%%%%%%%%%%%%%
\begin{eqnarray*}
 \mbox{where } \, \begin{cases}
\displaystyle T_{n,1}(q) = b_n(q)\sum_{\gamma\in \{-1,1\}, \, v\in
\mathcal{A}^{n}} \prod_{k=0}^{n-1} \ \mathbb{E}\big
(W_{v_{k+1}}(v|k)^{q+\gamma\eta_n}\big ),\\
\displaystyle T_{n,2}(q) =b_n(q) \sum_{\substack{\gamma\in \{-1,1\},\,
v,w\in \mathcal{A}^{n},\\ \delta(v,w)=1}}
\prod_{k=0}^{n-1}\mathbb{E}\big
(W_{v_{k+1}}(v|k)^{q}W_{w_{k+1}}(w|k)^{\gamma\eta_n}\big ).
\end{cases}
\end{eqnarray*}
One immediately gets
\begin{eqnarray*}
T_{n,1}(q)= \sum_{\gamma\in\{-1,1\}}b^{n\big
(\tau_\mu(q)+\gamma\eta_n
(\tau_\mu'(q)-\gamma\varepsilon_n)-\tau_\mu(q+\gamma\eta_n)\big
)}=2b^{-n (\varepsilon_n\eta_n+O(\eta_n^2) )},
\end{eqnarray*}
where $O(\eta_n^2)$ is uniform over $q\in K$ if
$\Vert\widetilde\eta\Vert_\infty$ is small enough (the twice
continuous differentiability of $\tau_\mu$ has been used).

\smallskip

Let $g_k$ and $d_k$ respectively stand for the word consisting of
$k$ consecutive zeros and the word consisting of $k$ consecutive
$b-1$. The estimation of $T_{n,2}(q)$ is achieved by using the following
identity:
%%%%%%%%%%%%%%%%%%%%
\begin{equation}
\label{decoupage}
\bigcup_{m=0}^{n-1} \ \bigcup_{u\in {\mathcal A}^{m}} \ \bigcup_{r\in
\{0,\dots,b-2\}} \big \{(u.r.d_{n-1-m}, u.(r+1).g_{n-1-m})\big \},
\end{equation}
%%%%%%%%%%%%%%%%%%
One has $T_{n,2}(q)=\mathcal{T}_n(q,-1)+\mathcal{T}_n(q,1)$ where for
$\gamma\in\{-1,1\}$
%%%%%%%%%%%%%%%%%%%%%%%%%%%%%%%%%%%%%%%%%%%%%%%%%%%%
\begin{eqnarray*}
\mathcal{T}_n(q,\gamma)&=&b_n(q)\sum_{\substack{v,w\in
\mathcal{A}^{n}\\\delta(v,w)=1}} \prod_{k=0}^{n-1}\mathbb{E}\big
(W_{v_{k+1}}(v|k)^{q}W_{w_{k+1}}(w|k)^{\gamma\eta_n}\big )\\
&=& b_n(q) \sum_{m=0}^{n-1}\,
\sum_{u\in\mathcal{A}^{m}}\, \sum_{r=0}^{b-2} \Theta_{n-1-m}(r)
\prod_{k=0}^{m-1} \, \mathbb{E}\big
(W_{u_{k+1}}(u|k)^{q+\gamma\eta_n}\big )\\
&=&b_n(q) \sum_{m=0}^{n-1} b^{-m\tau_\mu(q+\gamma\eta_n)} \,
\sum_{r=0}^{b-2}\Theta_{n-1-m}(r),
\end{eqnarray*}
%%%%%%%%%%%%%%%%%%%%%%%%%%%%%%%%%%%%%%%%%%%%%%%%%%%%
and where $\Theta_m(r)$ is defined by
\begin{eqnarray*}
\Theta_m(r) & = & \mathbb{E}\big (W_r^qW_{r+1}^{\gamma\eta_n}\big )\big
(\mathbb{E}(W_{b-1}^q)\big )^m\big
(\mathbb{E}(W_{0}^{\gamma\eta_n})\big )^m\\
&& \hspace{15mm} +
(\mathbb{E}\big (W_r^{\gamma\eta_n}W_{r+1}^{q}\big )\big
(\mathbb{E}(W_0^q)\big )^m\big
(\mathbb{E}(W_{b-1}^{\gamma\eta_n})\big )^m.
\end{eqnarray*}
All the components of $W$ are positive almost surely. Thus, by
definition (\ref{zero}) of $\widetilde\tau_\mu (q)=\tau_\mu(q)$, there
is a constant $c_K\in (0,1)$ such that for all $q\in K$ one has $\max
\big (\mathbb{E}(W_0^q),\mathbb{E}(W_{b-1}^q)\big )\le
c_Kb^{-\tau_\mu(q)}$. Moreover, if $\Vert\widetilde\eta\Vert_\infty$
is small enough, $\max \big (\mathbb{E}(W_0^{\gamma\eta_n}),
\mathbb{E}(W_{b-1}^{\gamma\eta_n})\big )\le (c^{-1}_K+1)/{2}$ (this
maximum goes to 1 when $\|\eta\|_\infty \rightarrow 0$). This yields
(since $c_K(c_K^{-1}+1) = c_K+1$)
\begin{eqnarray*}
\Theta_m(r)& \le & \big (\mathbb{E}\big (W_r^qW_{r+1}^{\gamma\eta_n}\big
)+(\mathbb{E}\big (W_r^{\gamma\eta_n}W_{r+1}^{q}\big )\big )\big
((c_K+1)/{2}\big) ^{m}b^{-m\tau_\mu(q) } \\ & \le & C_K\,\big
((c_K+1)/{2}\big) ^{m}b^{-m\tau_\mu(q) }.
\end{eqnarray*}
Consequently we  get
\begin{eqnarray*}
\mathcal{T}_n(q,\gamma)&\le & C_Kb^{n (\tau_\mu(q)+\gamma\eta_n
(\tau_\mu'(q)-\gamma\varepsilon_n)
)}b^{-(n-1)\tau_\mu(q+\gamma\eta_n)}\\ && \hspace{15mm} \times
\sum_{m=0}^{n-1}\big((c_K+1)/{2}\big) ^{m}b^{m
(\tau_\mu(q+\gamma\eta_n)-\tau_\mu(q) )}\\ &\le& C_Kb^{-n
(\varepsilon_n\eta_n+O(\eta_n^2) )}
\sum_{m=0}^{n-1}\big((c_K+1)/{2}\big) ^{m}b^{m
(\tau_\mu(q+\gamma\eta_n)-\tau_\mu(q) )}.
\end{eqnarray*}
The function $\tau_\mu$ is continuously differentiable. Hence the sum
$\sum_{m=0}^{n-1}((c_K+1)/{2}) ^{m}$ $b^{m
(\tau_\mu(q+\gamma\eta_n)-\tau_\mu(q) )}$ is uniformly bounded over
$n\ge 0$ and $q\in K$ if $\Vert\widetilde\eta\Vert_\infty$ is small
enough. Finally, if $\Vert\widetilde\eta\Vert_\infty$ is small enough
we also have $ \mathcal{T}_n(q,\gamma)\le C_Kb^{-n
(\varepsilon_n\eta_n+O(\eta_n^2))}$.

\noindent Going back to $T_n(q)$, we get $ T_n(q)\le C_K nb^{-n
(\varepsilon_n\eta_n+O(\eta_n^2) )}$, $\forall q\in K$. This
shows~(\ref{eq1}).

\smallskip

{\bf $\bullet$ The case $H=G$:} The proof follows similar lines as for
$F_n(q)$. The only new point it requires is the boundedness of
$\sup_{q\in K}\mathbb{E}(Y_q^{-h})$ for some $h>0$. In fact, we shall
need the following stronger property in the proof of
Theorem~\ref{cascades_ren}.
%%%%%%%%%%%%%%%%%%%%%%%%%%%%%%%%%%%%%%%%%%%%%%%%
\begin{lemma}\label{supmoments} 1.
For every compact subinterval $K$ of $J$, there exists $h>0$ such that
$ \mathbb{E}\big (\sup_{q\in K}Y_q^{-h}\big )<\infty $.\\
2. Assume that there exists $A>1$ such that with probability one,
$A^{-1}\le W_i$ (resp. $W_i\le A$) for all $0\le i\le b-1$. Then, for
every compact subinterval $K$ of $J\cap\mathbb{R}_+$
(resp. $J\cap\mathbb{R}_-$) there exist two constants depending on
$K$, $C_K>0$ and $\gamma_K\in (0,1)$ such that for all $x>0$ small
enough,
$$
\mathbb{P}(\inf_{q\in K}Y_q\le x)\le \exp\big
(-C_Kx^{-\gamma_K/(1-\gamma_K)}\big ).
$$\end{lemma}
%%%%%%%%%%%%%%%%%%%%%%%%%%%%%%%%%%%%%%%%%%%%%%%%
\begin{proof}1. Fix $K$ a 
compact subinterval of $J\cap\mathbb{R}_+$
(resp. $J\cap\mathbb{R}_-$).  For $w\in \mathcal{A}^*$ one can define
$Z_K(w)=\inf_{q\in K}Y_q(w)$ ($Z_K(\emptyset)$ is denoted $Z_K$). We
learn from \cite{Biggins,B2} that this infimum is positive since $q\mapsto
Y_q(w)$ is almost surely positive and continuous.

Let us also define $W_K(w)=\inf_{q\in K,0\le i\le b-1}W_{q,i}(w)$
($W_K(\emptyset)$ is denoted $W_K$). Since we assumed that $J$
contains a neighborhood of 0, there exists $h>0$ such that the moment
of negative order $-h$ of this random variable $W_K(w)$ is finite.

Moreover, with probability one, $\forall \,q\in J$ one has $
Y_q=\sum_{i=0}^{b-1}W_{q,i}(\emptyset)Y_q(i)$, hence
\begin{equation}
\label{ineg}
Z_K\ge W_K\sum_{i=0}^{b-1}Z_K(i).
\end{equation}
By construction, the random variables $Z_K(i)$, $0\le i\le b-1$, are
i.i.d. with $Z_K$, and they are independent of the positive random
variable $W_K$. Consequently, the Laplace transform of $Z_K$, denoted
$L:t\ \ge 0\mapsto \mathbb{E}(\exp (-t Z_K)\big )$, satisfies the
inequality
\begin{equation}\label{Laplace}
L(t)\le \mathbb{E}\Big (\prod_{i=0}^{b-1}L(W_K t)\Big )\quad (t\ge
0).
\end{equation}
Then, since $\mathbb{E}(W_K^{-h})<\infty$, using the approach of
 \cite{Mol} to study the behavior at $\infty$ of Laplace transforms
 satisfying an inequality like (\ref{Laplace}) (see also \cite{B1} and
 \cite{L3}) one obtains $\mathbb{E}(Z_K^{-h})<\infty$.

\smallskip

\noindent
2. It is a simple consequence of the proof of Theorem 2.5 in
   \cite{L3} (see also the proof of Corollary 2.5 in \cite{HOWA}) and
   of the fact that in this case, the random variable $W_K$ in
   (\ref{Laplace}) is lower bounded by a positive constant.

%%%%%%%%%%%%%%%%%%%%%%%%%%%%%%%%%%%%%%%%%%
%%%%%%%%%%%%%%%%%%%%%%%%%%%%%%%%%%%%%%%%%%
\end{proof}
%%%%%%%%%%%%%%%%%%%%%%%%%%%%%%%%%%%%%%%%%%

%%%%%%%%%%%%%%%%%%%%%%%%%%%%%%%%%%%%%%%%%%%
%%%%%%%%%%%%%%%%%%%%%%%%%%%%%%%%%%%%%%%%%%%
\subsection{Proof of Theorem~\ref{cascades_ren}}
Fix $K$ a compact subinterval of $J$. The computations
performed to prove Theorem~\ref{cascades} yield
(\ref{crucial}). Thus there are two constants $C>0$ and $\beta>0$ as
well as a sequence $\widetilde\eta =(\eta_n)_{n\ge 1}\in
\mathbb{R}_+^{\mathbb{N}^*}$ such that for every $j,n\ge 1 $, $q\in K$
and $w\in \mathcal{A}^j$,
\begin{equation}\label{good}
\mathbb{E}\Big (\sup_{q\in K}S^{N,\varepsilon_n,\eta_n}\big
(\mu_q^{w},\mu^{w},\tau_\mu'(q)\big )\Big )
\le C n b^{-\beta n\varepsilon^2_n}.
\end{equation}

In order to apply Proposition~\ref{ren}, let us define

\smallskip
\noindent
$\bullet$ $\Lambda= K$, $\{(m_\lambda
^w,\mu_\lambda^w)\}_{w\in\mathcal{A}^*,\lambda\in \Lambda}=\{\mu_q
^w,\mu^w\}_{w\in\mathcal{A}^*,q\in K}$ and
$\{\alpha_\lambda\}_{\lambda\in \Lambda}=\{\tau_\mu'(q)\}_{q\in K}$,

\smallskip
\noindent
$\bullet$ For $w\in \mathcal{A}^*$ and $n\ge 1$, 
\begin{equation}
\label{defv}
U^w=\inf_{q\in K}\Vert\mu_q^w\Vert\quad\mbox{and}\quad
V^w_n=\sup_{q\in K}S^{N,\varepsilon_n,\eta_n}\big
(\mu_q^{w},\mu^{w},\tau_\mu'(q)\big ),
\end{equation}

\noindent
$\bullet$ For every $j\ge 1$, $\psi_j(\widetilde\eta)=1$ and
$\rho_j=\log (j)^{1+\eta}$.

\smallskip
\noindent
$\bullet$ Fix $\eta>0$ and $\eta'>2\eta$.  For every $j\ge 1$, we set
$\mathcal{S}_j=\big [\exp \big (\big (j\log
(j)^{\eta}\big)^{\frac{1}{1+2\eta}}\big ) \big ]$ if $ {\log
(j)^{-\eta}} \geq \varepsilon_j\geq j^{-1/2}\log (j)^{1/2+\eta}$ and
$\mathcal{S}_j=\big [j\log (j)^{\eta'}\big ]$ if $\varepsilon_j\geq
{\log (j)^{-\eta}}$.

\smallskip

Now, on the one hand, Lemma~\ref{supmoments}.2 implies that
\begin{equation}\label{good2}
u_j:=b^j\mathbb{P}(U^w\le b^{-\rho_j})\le b^j\exp \big (-C_K
b^{\frac{\gamma_K}{1-\gamma_K}(\log j)^{(1+\eta)}}\big ).
\end{equation}
Moreover, $\sum_{j\ge 1}u_j<\infty$.  On the other hand, for some
$\alpha>0$, for any $w\in\mathcal{A}^*$ one has
$$
v_j:=\sum_{n\ge \mathcal{S}_j}\mathbb{E}(V^w_n)\le \sum_{n\ge
\mathcal{S}_j} C n  b^{-\beta
n\varepsilon^2_n}=O\big (b^{-j\log (j)^\alpha}\big ) .
$$ The sequence $\rho_j$ has been chosen so that $\sum_{j\ge
1}b^jb^{\rho_j}v_j<\infty$.  Consequently, Proposition~\ref{ren}
yields the desired upper bound for the growth speed
$GS\big(\mu^w_q,\mu^w,\tau_\mu'(q), N,\widetilde\varepsilon\big)$.

Changing the measures $\{(m_q ^w,\mu_q^w)\}_{w\in\mathcal{A}^*,q\in
K}$ into $\{\mu_q^w,\mu_q^w\}_{w\in\mathcal{A}^*,q\in K}$ and the
exponents $\{\tau_\mu'(q)\}_{q\in K}$ into
$\{\tau_\mu'(q)q-\tau_\mu(q)\}_{q\in K}$, the same arguments yield the
conclusion for $GS\big(\mu^w_q,\mu^w_q,\tau_\mu'(q)q-\tau_\mu(q)),
N,\widetilde\varepsilon\big)$.

%%%%%%%%%%%%%%%%%%%%%%%%%%%%%%%%%
%%%%%%%%%%%%%%%%%%%%%%%%%%%%%%%%%%%%%%%%%%%
\subsection{Proof of Theorem~\ref{cascades_ren2}} 

We only prove the results for the control of\\
$GS\big(\mu^w_q,\mu^w,\tau_\mu'(q), N,\widetilde\varepsilon\big)$ by
$\mathcal{S}_j$, since
$GS\big(\mu^w_q,\mu^w_q,\tau_\mu'(q)q-\tau_\mu(q)),
N,\widetilde\varepsilon\big)$ is controlled by using the same
approach.

\smallskip

{\bf (1)} Recall that $(\Omega,\mathcal{B},\mathbb{P})$ denotes the
probability space on which the random variables in this
Section~\ref{CASCADES} are defined. Let us consider on~$\mathcal
{B}\otimes \mathcal{B}([0,1])$ the so-called Peyri\`ere probability
$\mathcal{Q}_q$ \cite{KP}
$$
\mathcal{Q}_q(A)=\mathbb{E}\Big (\int_{[0,1]}\mathbf{1}_A(\omega,t)\,
\mu_q(t) \Big )\quad (A\in \mathcal {B}\otimes \mathcal{B}([0,1])).
$$
It is important to notice that by construction $\mathcal{Q}_q$-almost
surely means $\mathbb{P}$-almost surely, $\mu_q(\omega)$-almost
everywhere.

Fix $\widetilde\eta$ as in the proof of
Theorem~\ref{cascades_ren}. Also, for $j\ge 1$ let $\rho_j=\log
(j)^{1+\eta}$, and let $\mathcal{S}_j=\big [j\log (j)^{-\kappa}\big
]$.  Now, for $j\ge 0$ and $n\ge 1$ define on $\Omega\times [0,1)$ the
random variables
\begin{eqnarray*}
\mathbf{U}^{(j)}(\omega,t)&=&\Vert\mu_q^{w^{(j)}(t)}(\omega)\Vert.\\
\mbox{and }
\mathbf{V}^{(j)}_n(\omega,t)&=&S_n^{N,\varepsilon_n,\eta_n} \big
(\mu_q^{w^{(j)}(t)}(\omega),\mu^{w^{(j)}(t)}(\omega),\tau_\mu'(q)\big ).
\end{eqnarray*}
We can use the proof of Proposition~\ref{ren} to claim that it is
enough to prove that
%%%%%%%%%%%%%
$$\exists \ h\in (0,1],\  \sum_{j\ge 0}\mathcal{Q}_q\big
(\mathbf{U}^{(j)}\le b^{-\rho_j}\big )<\infty \mbox{ and } \sum_{j\ge
0}b^{\rho_jh}\mathbb{E}_{\mathcal{Q}_q}\big (\big (\sum_{n\ge
\mathcal{S}_j}\mathbf{V}^{(j)}_n\big )^h\big )<\infty.
$$
%%%%%%%%%%%%%%%%%%%%%%%
The main difference with the proofs of Proposition~\ref{ren} and
Theorem~\ref{cascades_ren} is that here we do not seek for a result
valid uniformly over the $w$ of the same generation $j$, but only for
a result valid for $w^{(j)}(t)$, for $\mu$-almost every $t$. As a
consequence we must control only one pair of random variables
$(\mathbf{U}^{(j)},\mathbf{V}^{(j)})$ on each generation instead of
$b^j$. This allows to slow down $\mathcal{S}_j$.

Fix $h\in (0,1)$. Since $x\mapsto x^h$ is sub-additive on
$\mathbb{R}_+$, one has
%%%%%%%%%%%%%%%%
$$
\mathbb{E}_{\mathcal{Q}_q}\Big (\Big (\sum_{n\ge
\mathcal{S}_j}\mathbf{V}^{(j)}_n\Big )^h\Big )\le \sum_{n\ge
\mathcal{S}_j}\mathbb{E}_{\mathcal{Q}_q}\Big (\big
(\mathbf{V}^{(j)}_n\big )^h\Big ).
$$
%%%%%%%%%%%%%%%%

%Let us estimate $\mathbb{E}_{\mathcal{Q}_q}\Big (\big
%(\mathbf{V}^{(j)}_n\big )^h\Big )$.  
For $\omega\in\Omega^*$, $j\ge 1$ and $n\ge 1$, by definition of the
measures $\mu_q$ and $\mu_q^w$, and since $\big
(\mu_q^{w^{(j)}(t)}(\omega),\mu^{w^{(j)}(t)}(\omega)\big )$ does not depend on
$t\in I_w$, one has
%%%%%%%%%%%%%%%%
\begin{eqnarray*}
&&\int_{[0,1]}\big (\mathbf V^{(j)}_n(\omega,t)\big )^h\mu_q(\omega)(dt)
=\sum_{w\in\mathcal{A}^j}\int_{I_w}\big (\mathbf V^{(j)}_n(\omega,t)\big
)^h\mu_q(\omega)(dt)\\
&=& \sum_{w\in\mathcal{A}^j}\prod_{k=0}^{j-1}W_{q,w_{k+1}}(
w|k)\int_{I_w}\big (\mathbf V^{(j)}_n(\omega,t)\big )^h\mu^w_q(\omega)\circ
f_{I_w}^{-1}(dt)\\
&=& \sum_{w\in\mathcal{A}^j}\Big (\prod_{k=0}^{j-1}W_{q,w_{k+1}}(
w|k)\Big ) \big (V_n^w\big )^h\, \Vert\mu_q^w\Vert,
\end{eqnarray*}
%%%%%%%%%%%%%%%%
where $V_n^w= S_n^{N,\varepsilon_n,\eta_n} \big
(\mu_q^{w}(\omega),\mu^{w}(\omega),\tau_\mu'(q)\big )$ is defined as
in the proof of Theorem~\ref{cascades_ren}. The above sum is a random
variable on $(\Omega,\mathcal{B},\mathbb{P})$. In addition, in each of
its terms, the product is independent of $\big (V_n^w\big )^h\,
\Vert\mu_q^w\Vert$. Moreover, the probability distribution of $\big
(V_n^w\big )^h\, \Vert\mu_q^w\Vert$ does not depend on
$w$. Consequently, using the martingale property of
$\Vert\mu_{q,n}\Vert$, one gets
%%%%%%%%%%%%%%%%
$$
\mathbb{E}_{\mathcal{Q}_q} \left (\big (\mathbf{V}^{(j)}_n\big
)^h\right )=\mathbb{E}\left (\big (V_n^w\big )^h\,
\Vert\mu_q^w\Vert\right ),
$$
%%%%%%%%%%%%%%%%
where $w\in \mathcal{A}^j$. Let $p=1/(1-h)$. The H\"older inequality
yields
%%%%%%%%%%%%%%%%
$$
\mathbb{E} \left(\big (V_n^w\big )^h\, \Vert\mu_q^w\Vert\right )\le
\big  (\mathbb{E}\big (V_n^w\big )\big )^{h}\, \mathbb{E}\big
(\Vert\mu_q\Vert^p\big )^{1/p}.
$$
%%%%%%%%%%%%%%%% 
Finally, $p$ is fixed close enough to 1 so that $\mathbb{E}\big
(\Vert\mu_q\Vert^p\big )<\infty$ (see the proof of
Lemma~\ref{Momentspos} for the existence of such a $p$).  Then
(\ref{good}) yields $\sum_{j\ge
1}b^{\rho_jh}\sum_{n\ge\mathcal{S}_j}\big (\mathbb{E}\big (V_n^w\big
)\big )^{h}<\infty$, hence the conclusion.

\medskip

Similar computations as above show that for every $j\ge 1$,
%%%%%%%%%%%%%%%%
$$
\mathcal{Q}_q\big (\mathbf{U}^{(j)}\le b^{-\rho_j}\big
)=\mathbb{E}\big (\mathbf{1}_{\{Y_q\le b^{-\rho_j}\}}Y_q\big )\le
b^{-\rho_j}\mathbb{P}(Y_q\le b^{-\rho_j}).
$$
%%%%%%%%%%%%%%%%
It follows from Lemma~\ref{supmoments}.1 that for some $h>0$ one has
$\mathbb{P}(Y_q\le x)=O(x^h)$ as $x\to 0$. This implies $\sum_{j\ge
1}\mathcal{Q}_q\big (\mathbf{U}^{(j)}\le b^{-\rho_j}\big )<\infty$.
\medskip

\noindent
{\bf (2)} The proof is similar to the one of (1). It is enough to prove the
    result for a compact subinterval $K$ of $J$ instead of $J$. Fix
    such an interval $K$. The idea is now to consider on $\big
    (K\times \Omega\times [0,1],
    \mathcal{B}(K)\otimes\mathcal{B}\otimes\mathcal{B}([0,1])\big )$
    the probability distribution $\mathcal{Q}_K$
%%%%%%%%%%%%%%%%
$$
\mathcal{Q}_K(A)=\int_K\left (\mathbb{E}_{\mathcal{Q}_q}
\mathbf{1}_{A}(q,\omega,t)\right )\, \frac{dq}{|K|}\quad \big (A\in
\mathcal{B}(K)\otimes\mathcal{B}\otimes\mathcal{B}([0,1])\big ).
$$
%%%%%%%%%%%%%%%%
 Then $\mathbf{U}^{(j)}(q,\omega,t)$ and
$\mathbf{V}^{(j)}_n(q,\omega,t)$ are redefined as
$\displaystyle\mathbf{U}^{(j)}(q,\omega,t) =
\Vert\mu_q^{w^{(j)}(t)}(\omega)\Vert$ and $\displaystyle
\mathbf{V}^{(j)}_n(q,\omega,t) = S_n^{N,\varepsilon_n,\eta_n} \big
(\mu_q^{w^{(j)}(t)}(\omega),\mu^{w^{(j)}(t)}(\omega),\tau_\mu'(q)\big
)$.

%%%%%%%%%%%%%%%%
Since there exists $p>1$ such that $M=\sup_{q\in K}\mathbb{E}\big
(\Vert\mu_q\Vert^p\big )^{1/p}<\infty$ (see the proof of
Lemma~\ref{Momentspos}), the computations performed above yield
%%%%%%%%%%%%%%%%
\begin{eqnarray*}
\sum_{j\ge 0}b^{\rho_j h}\sum_{n\ge
\mathcal{S}_j}\mathbb{E}_{\mathcal{Q}_K} \left(\big
(\mathbf{V}^{(j)}_n\big )^h\right )\le |K|\,M\sum_{j\ge
1}b^{\rho_jh}\sum_{n\ge\mathcal{S}_j}\Big (\mathbb{E}\big (V_n^w\big
)\Big )^{h}<\infty.
\end{eqnarray*}
%%%%%%%%%%%%%%%%
Finally, $ \sum_{j\ge 0}\mathcal{Q}_K\big (\mathbf{U}^{(j)}\le
b^{-\rho_j}\big )\le |K|\sum_{j\ge 1}b^{-\rho_j}\mathbb{P}\big
(\inf_{q\in K}Y_q\le b^{-\rho_j}\big )$, which is finite by item 1. of
Lemma \ref{supmoments}.

%%%%%%%%%%%%%%%%%%%%%%%%%%%%%%%%%%%%%%%%%%%
%%%%%%%%%%%%%%%%%%%%%%%%%%%%%%%%%%%%%%%%%%%
\subsection{Proof of Theorem~\ref{cascades_dev}}

 We assume without loss of generality that $K$ contains the point
 1. Define $q_K=\max \{|q|: q\in K\}$.  Recall that for $j\ge 0$ and
 $n\ge 1$, if $(w,v)\in\mathcal{A}^j\times\mathcal{A}^n$ and $q\in K$
 then
$$
\mu^w(I_v)^q=\mu^w_q(I_v) b^{-n\widetilde\tau_\mu(q)}\frac{Y(wv)^q}{Y_q(wv)}.
$$
As a consequence,
\begin{eqnarray}\label{GIBBS1}
Y_q(w)b^{-n\widetilde\tau_\mu(q)}\inf_{q\in K,v\in
\mathcal{A}^n}\frac{Y(wv)^q}{Y_q(wv)}\le b^{-n\tau^w_n(q)}\\
\mbox{and } \  b^{-n\tau^w_n(q)}\le
Y_q(w)b^{-n\widetilde\tau_\mu(q)}\sup_{q\in K,v\in
\mathcal{A}^n}\frac{Y(wv)^q}{Y_q(wv)}.
\end{eqnarray}
Let us fix $\delta\in (0,1)$ and $\theta>0$ such that the conclusions
of Propositions \ref{cont1} and \ref{cont2} below hold. Then, with
probability one, for $j$ large enough, for every $ w\in
\mathcal{A}^j$, $q\in K$ and $n\ge j^\delta$, one has $
b^{-n\tau^w_n(q)}\le
Y_q(w)b^{-n\widetilde\tau_\mu(q)}n^{(q_K+1)\theta}.  $ This yields
%%%%%%%%%%%%%%%%
\begin{eqnarray*}
\mathcal N_n(\mu^w,\widetilde\tau'_\mu(q),\varepsilon_n) & \le  & b^{n\big
(\widetilde\tau'_\mu(q)q-\tau^w_n(q)+sgn(q)q\varepsilon_n\big )}\\
& \le & Y_q(w)
b^{n\big
(\widetilde\tau'_\mu(q)q-\widetilde\tau_\mu(q)+sgn(q)q\varepsilon_n\big
)}n^{(q_K+1)\theta}.
\end{eqnarray*}
%%%%%%%%%%%%%%%%

On the other hand, due to Theorem~\ref{cascades_ren} and
Proposition~\ref{cont1}, there exists $\theta>0$ such that, with
probability one, for $j$ large enough, for all $w\in \mathcal{A}^*$
and $q\in K$
%%%%%%%%%%%%%%%%
$$ \mu_q^{w}\Big (E^{\mu^{w}}_{\widetilde\tau'_\mu(q),\mathcal{S}_j}
(0,\widetilde\varepsilon)\bigcap E^{\mu_q^{w}}_{\widetilde\tau'_\mu(q)
q- \widetilde\tau_\mu(q), \mathcal{S}_j}(0,\widetilde\varepsilon)\Big
)\ge \Vert\mu_q^{w}\Vert/2,
$$ which equals ${Y_q(w)}/{2}$. Thus $b^{n(\widetilde\tau'_\mu(q)q-
\widetilde\tau_\mu(q)-\varepsilon_n) }{Y_q(w)}/{2}\le \mathcal
N_n(\mu^w,\widetilde\tau'_\mu(q),\varepsilon_n)$ for every $n\ge
\mathcal{S}_j$. Moreover, for $j$ large enough one has $j^\delta \le
\mathcal{S}_j$, and then for $n$ large enough $\sup_{q\in
K}sgn(q)q\varepsilon_n+(q_K+1)\theta\, {\log_b(n)}/{n}$ is controlled
by $(1+q_K)\varepsilon_n$. The conclusion follows.
%%%%%%%%%%%%%%%%%%%%%%%%%%%%%%%%%%%%%%%%%%%
%%%%%%%%%%%%%%%%%%%%%%%%%%%%%%%%%%%%%%%%%%%
%%%%%%%%%%%%%%%%%%%%%%%%%%%%%%%%%%%%%%%%%%%
%%%%%%%%%%%%%%%%%%%%%%%%%%%%%%%%%%%%%%%%%%%
\subsection{Proof of Theorem \ref{convergence}}

Let us begin with three technical lemmas.
%%%%%%%%%%%%%%%%%%%%%%%%%%%%%%%%%%%%%%%%%%
\begin{lemma}\label{Momentspos}
Let us assume that $J=\mathbb{R}$. For every compact subinterval $K$
of $\mathbb{R}$, there exist $C_K,c_K>0$ such that
$$ \mbox{for every }x\ge 1, \ \sup_{q\in K}\mathbb{P}(Y_q\ge x)\le
C_K\exp\big (-c_Kx).
$$
\end{lemma}
%%%%%%%%%%%%%%%%%%%%%%%%%%%%%%%%%%%%%%%%%%

%%%%%%%%%%%%%%%%%%%%%%%%%%%%%%%%%%%%%%%%%%

%%%%%%%%%%%%%%%%%%%%%%%%%%%%%%%%%%%%%%%%%%

%%%%%%%%%%%%%%%%%%%%%%%%%%%%%%%%%%%%%%%%%%
\begin{proof}
 Let us begin with the following properties involved in the proofs of
several statements of Section~\ref{CASCADES}: it is known (see
\cite{KP,DL}) that if $q\in {J}$ and $h>1$ one has
$\mathbb{E}(Y_q^h)<\infty$ if and only if $\mathbb{E}\Big
(\sum_{i=0}^{b-1}W_{q,i}^h\Big )<1$, that is $\widetilde\tau_\mu
(qh)-h\widetilde\tau_\mu(q)>0$. Moreover, one deduces from the proofs of
Theorem III.B. and Theorem VI.A.b. of \cite{B1}: for every compact
subinterval $K$ of $J$, there exists $h>1$ such that $\sup_{q\in
K}\mathbb{E}(Y_q^h)<\infty$.

The property $J=\mathbb{R}$ is equivalent to the fact that the mapping
$q\mapsto {\widetilde \tau_\mu(q)}/{q}$ is increasing on
$\mathbb{R}_-^*$ and $\mathbb{R}_+^*$. As a consequence, one has
$\widetilde\tau_\mu (qh)-h\widetilde\tau_\mu(q)>0$ for all $q\in\mathbb{R}$
and $h>1$, that is $\mathbb{E}(Y_q^h)<\infty$. Also, one has $\Vert
W_{q,i}\Vert_\infty \le 1$ for all $q\in\mathbb{R}$ and $0\le i\le
b-1$.

Let us now fix $K$ a compact subset of $\mathbb{R}$. Let us then
introduce $\displaystyle t_k(q)={\mathbb{E}(Y_q^k)}/{k !}$ for
$q\in K$ and $k\ge 1$, and $t_0(q)=1$. The proof of Theorem 4.1 (a) in
\cite{L1} yields that for every $k\geq 2$, for every $q\in K$, 
$$ t_k(q)\le \ c_K\sum_{(k_0,\dots,k_{b-1}): \,0\le k_i\le
k-1\mbox{ {\tiny and }} k_0+\dots +k_{b-1}=k} \ \ \prod_{i=0}^{b-1}t_{k_i}(q),
$$
where $c_K=\sup_{q\in K}\sup_{k\ge 2}\big
(1-b^{-\widetilde\tau_\mu(kq)+k\widetilde\tau_\mu (q)}\big )^{-1}$. One sees
that $c_K=\sup_{q\in
K}\big(1-b^{-\widetilde\tau_\mu(2q)+2\widetilde\tau_\mu (q)}\big
)^{-1}<\infty$. Hence, if $\widetilde t_k= \sup_{q\in
  K}t_k(q)$, one has 
%%%%%%%%%%%%%%%%%%%%%%%%%%%%%%%%%%%%%%%%%%
\begin{equation*}
\forall \, k\ge 2,\ \ \widetilde t_k\le \
c_K\sum_{{(k_0,\dots,k_{b-1}): \,0\le k_i\le k-1\mbox{ {\tiny and }}
k_0+\dots +k_{b-1}=k}} \ \ \prod_{i=0}^{b-1}\widetilde t_{k_i} .
\end{equation*}
%%%%%%%%%%%%%%%%%%%%%%%%%%%%%%%%%%%%%%%%%%
Since $\widetilde t_0=\widetilde t_1=1$, Lemma 2.6 of \cite{GMW}
yields $\limsup_{k\to\infty}\widetilde t_k^{\frac{1}{k}}<\infty
$. This implies the existence of a constant $C>0$ such that
%%%%%%%%%%%%%%%%%%%%%%%%%%%%%%%%%%%%%%%%%%
\begin{equation*}
 \forall \, k\ge 1, \ \ \sup_{q\in K}\mathbb{E}(Y_q^k)\le C^k k!.
\end{equation*}
%%%%%%%%%%%%%%%%%%%%%%%%%%%%%%%%%%%%%%%%%%
Now, fix $c_K\in (0,C^{-1})$. For $x>0$ one has 
%%%%%%%%%%%%%%%%%%
\begin{eqnarray*}
\sup_{q\in K}\mathbb{P}(Y_q\ge x)\le e^{-c_Kx} \sup_{q\in
K}\mathbb{E}\big (e^{c_KY_q}\big )&\le &e^{-c_Kx}\sum_{k=0}^\infty
c_K^k {\sup_{q\in K}\mathbb{E}(Y_q^k)}/{k!}\\
& \le &
(1-c_KC)^{-1}e^{-c_Kx}.
\end{eqnarray*}
\end{proof}
%%%%%%%%%%%%%%%%%%%%%%%%%%%%%%%%%%%%%%%%%%

\begin{remark}
We are not able to control $\mathbb{P}(\sup_{q\in K}Y_q\ge x)$ at
$\infty$. This is the reason why next Lemmas \ref{Holder} and
\ref{pp1} are needed to prove Proposition~\ref{cont2}.
\end{remark}

For $n\ge 1$ let $Q_n$ be the set of dyadic numbers of generation $n$.
%%%%%%%%%%%%%%%%%%%%%%%%%%%%%%%%%%%%%%%%%%
\begin{lemma}\label{Holder}
Let $K$ be a compact subinterval of $J$. Let $\eta >0$. There exists
  $\alpha\in (0,1)$ and $\delta\in (0,1)$ such that, with probability
  one,\\
{\bf 1.} for $j$ large enough, $\forall\ w\in \mathcal{A}^j$, $\forall\
n\ge [j^{1+\eta}]$, $\forall\ q,q'\in Q_n$ such that $|q-q'|= 2^{-n}$,
one has $\big |Y_q^w-Y_{q'}^w\big |\le |q'-q|^\alpha$.\\
{\bf 2.} for $j$ large enough, $\forall n\ge j^{\delta}$, $\forall \ w\in
\mathcal{A}^j,\ \forall v\in \mathcal{A}^n$, $\forall \ m\ge
[n^{1+\eta}]$, for all $q,q'\in Q_m$ such that $|q'-q|= 2^{-m} $, one
has $\big |Y_q^{wv}-Y_{q'}^{wv}\big |\le |q'-q|^{\alpha}$.
\end{lemma}
%%%%%%%%%%%%%%%%%%%%%%%%%%%%%%%%%%%%%%%%
\begin{proof}
By Theorem VI.A.b. $i)$ of \cite{B1}, $\exists \ h>1, \ C_K>0$ such that 
\begin{equation}\label{holder1}
\mbox{for all }(q,q')\in K^2,\ \mathbb{E}\left (|Y_q-Y_{q'}|^h\right)\le
C_K |q-q'|^h.
\end{equation}
For $n\ge 1$, let $\widetilde Q_n$ be the set of pairs $(q,q')\in Q_n$
such that $|q-q'|=2^{-n}$, and let $\alpha\in
(0,(h-1)/h)$. Using (\ref{holder1}) and the Markov
inequality, one gets
\begin{eqnarray*}
p_n&:=&\mathbb{P}\left (\exists\ (q,q')\in \widetilde Q_n,\
|Y_q-Y_{q'}|\ge |q-q'|^\alpha\right )\\
&\le &\sum_{(q,q')\in \widetilde Q_n}\mathbb{P}\Big (|Y_q-Y_{q'}|\ge
|q-q'|^\alpha\Big ) \leq 2|K|2^n C_K2^{n\alpha h}2^{-nh}.
\end{eqnarray*}
Let us fix $\eta>0$ and $\delta\in (({1+\eta})^{-1},1)$.  $\sum_{j\ge
1}b^j\sum_{n\ge [j^{1+\eta}]}p_n<\infty$ implies item {\bf 1.} of
Lemma \ref{Holder} by the Borel-Cantelli lemma. Also, item {\bf 2.}
follows from the fact that $\sum_{j\ge 1}b^j\sum_{n\ge
j^\delta}b^n\sum_{m\ge [n^{1+\eta}]}p_m<\infty$.
\end{proof}
%%%%%%%%%%%%%%%%%%%%%%%%%%%%%%%%%%%%%%%%
%%%%%%%%%%%%%%%%%%%%%%%%%%%%%%%%%%%%%%%%

%%%%%%%%%%%%%%%%%%%%%%%%%%%%%%%%%%%%%%%%%%%
\begin{lemma}\label{pp1} Under the assumptions of
  Theorem~\ref{cascades_dev}, let $K$ be a compact subinterval of
$J=\mathbb{R}$. Let $\eta>0$. There exist $\delta\in (0,1)$ and
$\theta >1$ such that, with probability one, for $j$ large enough,

\noindent
1. $\forall\, w\in \mathcal{A}^j$,
$\sup_{q\in Q_{[j^{1+\eta}]}\cap K}Y^w_q\le j^\theta$.\\
2. $\forall\, n\ge j^\delta$, $\forall \, (v,w)\in 
\mathcal{A}^n \times \mathcal{A}^j$, $\sup_{q\in Q_{[n^{1+\eta}]}\cap
K}Y^{wv}_q\le~n^\theta$.
\end{lemma}
%%%%%%%%%%%%%%%%%%%%%%%%%%%%%%%%%%%%%%%%%%
\begin{proof}
Fix $\theta>1+\eta$. For $q\in K$ and $j\ge 1$ define
$p_j(q)=\mathbb{P}(Y_q\ge j^\theta)$. By Lemma~\ref{Momentspos},
\begin{eqnarray*} 
\forall \, j\ge 1,
 \ \ \mathbb{P}\Big(\sup_{q\in Q_{[j^{1+\eta}]}\cap K}Y_q\ge j^ \theta
\Big) & \le & \sum_{q\in Q_{[j^{1+\eta}]}\cap K}p_j(q) \\ & \le&
p_j:=2C_K|K|2^{j^{1+\eta}}\exp\big (-c_Kj^{\theta_K}).
\end{eqnarray*} 
We let the reader verify that $\sum_{j\ge 1}b^jp_j<\infty$ and
$\sum_{j\ge 1}b^j\sum_{n\ge j^\delta}b^np_n<\infty$ if $\delta\in
(\theta_K^{-1},1)$. This yields items 1. and 2. of Lemma
\ref{pp1}.
\end{proof}
%%%%%%%%%%%%%%%%%%%%%%%%%%%%%%%%%%%%%%%%%%

Let us now finish the proof of Theorem \ref{convergence}. It is a
   consequence of (\ref{GIBBS1}) and of the next Propositions
   \ref{cont1} and \ref{cont2}.

%%%%%%%%%%%%%%%%%%%%%%%%%%%%%%%%%%%%%%%%%%%
\begin{proposition}\label{cont1}
Under the assumptions of Theorem~\ref{cascades_ren}, let $K$ be a
compact subinterval of $J\cap\mathbb{R}_+$
(resp. $J\cap\mathbb{R}_-$). There exist $\theta >0$ and $\delta\in
(0,1)$ such that, with probability one, for $j$ large enough

\smallskip
\noindent
1. $\forall\ w\in \mathcal{A}^j$, one has
$\inf_{q\in K}Y^w_q\ge j^{-\theta}$.\\
2.  $\forall n\ge j^{\delta}$, $\forall \ w\in \mathcal{A}^j,\ \forall
v\in \mathcal{A}^n$, one has $\inf_{q\in K}Y^{wv}_q\ge~n^{-\theta}$.
\end{proposition}
%%%%%%%%%%%%%%%%%%%%%%%%%%%%%%%%%%%%%%%%%%%
\begin{proof} Fix $\theta>1$ such that $\theta_K=\theta
\frac{\gamma_K}{1-\gamma_K}>1$, where $\gamma_K$ is as in
Lemma~\ref{supmoments}. Let us also define
$p_j=\mathbb{P}\left(\inf_{q\in K}Y_q<j^{-\theta}\right)$.

We let the reader verify, using Lemma~\ref{supmoments}, that
$\sum_{j\ge 1}b^jp_j<\infty$ and if $\delta\in (\theta_K^{-1},1)$ then
$\sum_{j\ge 1}b^j\sum_{n\ge j^\delta}b^np_n<\infty$. This yields
1. and 2.
\end{proof}
%%%%%%%%%%%%%%%%%%%%%%%%%%%%%%%%%%%%%%%%
\begin{proposition}\label{cont2} Under the assumptions of
  Theorem~\ref{cascades_dev}, let $K$ be
a compact subinterval of $J$. There exist $\theta >0$ and $\delta\in
(0,1)$ such that, with probability one, for $j$ large enough

\smallskip
\noindent
1. $\forall\ w\in \mathcal{A}^j$, one has
$\sup_{q\in K}Y^w_q\le j^\theta$.\\
2. $\forall n\ge j^{\delta}$, $\forall \ w\in
\mathcal{A}^j,\ \forall v\in \mathcal{A}^n$, one has $\sup_{q\in
K}Y^{wv}_q\le n^\theta$.
\end{proposition}
%%%%%%%%%%%%%%%%%%%%%%%%%%%%%%%%%%%%%%%%%%

%%%%%%%%%%%%%%%%%%%%%%%%%%%%%%%%%%%%%%%%%%
\begin{proof}
We assume without loss of generality that the end points of $K$ are
dyadic numbers.

It is standard (see the proof of Kolmogorov theorem in \cite{KARSH})
that Lemma~\ref{Holder} implies that there exists a constant $C_K>0$
such that, with probability one\\
1. for $j$ large enough, $\forall\ w\in \mathcal{A}^j$, $\forall\
q,q'\in K$ such that $|q-q'|\le 2^{-[j^{1+\eta}]}$, one has $\big
|Y_q^w-Y_{q'}^w\big |\le C_K|q'-q|^\alpha$.\\
2. for $j$ large enough, $\forall n\ge j^{\delta}$, $\forall \ (v,w)\in
\mathcal{A}^n \times \mathcal{A}^j$, for all $q,q'\in K$ such
that $|q'-q| \le 2^{-[n^{1+\eta}]}$, one has $\big
|Y_q^{wv}-Y_{q'}^{wv}\big |\le C_K |q'-q|^{\alpha}$.

Then, Lemma~\ref{pp1} concludes the proof.
\end{proof} 
%%%%%%%%%%%%%%%%%%%%%%%%%%%%%%%%%%%%%%%%%%

%%%%%%%%%%%%%%%%%%%%%%%%%%%%%%%%%%%%%%%%%%%%%%
%%%%%%%%%%%%%%%%%%%%%%%%%%%%%%%%%%%%%%%%%%%%%
%%%%%%%%%%%%%%%%%%%%%%%%%%%%%%%%%%%%%%%%%%%%%%
%%%%%%%%%%%%%%%%%%%%%%%%%%%%%%%%%%%%%%%%%%%%%
%%%%%%%%%%%%%%%%%%%%%%%%%%%%%%%%%%%%%%%%%%%%%%
%%%%%%%%%%%%%%%%%%%%%%%%%%%%%%%%%%%%%%%%%%%%%
\section{The version of Theorem~\ref{cascades_ren2} needed to get (\ref{DIM})}
\label{ubi}

Now, let $\{(x_n,\lambda_n)\}_{n\ge 1}$ be a sequence in $[0,1]\times
(0,1]$ such that $\lim_{n\to\infty}\lambda_n=0$.

For every $t\in (0,1)$, $k\ge 1$ and $r\in (0,1)$ let us define
$$
\mathcal{B}_{k,r}(t)=\big \{B(x_n,\lambda_n):\ t\in B(x_n,r\lambda_n),\ 
\lambda_n\in (b^{-(k+1)},b^{-k}]\big \}
$$ ($B(y,r)$ denotes the closed interval centered at $y$ with radius
$r$). Notice that this set may be empty. Then, if $\xi>1$ and
$B(x_n,\lambda_n)\in \mathcal{B}_{k,1/2}(t)$, let
$\mathcal{B}_k^\xi(t)$ be the set of $b$-adic intervals of maximal
length included in $B(x_n,\lambda^\xi_n)$.
%%%%%%%%%%%%%%%%%%%%%%%%%%%%%%%%%%%%%%%%%%%%%
%%%%%%%%%%%%%%%%%%%%%%%%%%%%%%%%%%%%%%%%%%%%%
The next result is key to build a generalized Cantor set of Hausdorff
dimension $\ge \tau_\mu^*(\alpha)/\xi$ in the set $K(\alpha,\xi)$ (\ref{K}).

%%%%%%%%%%%%%%%%%%%%%%%%%%%%%%%%%%%%%%%%%%%%%
\begin{theorem}\label{analogue2}
Suppose that $\limsup_{n\to\infty}B(x_n,\lambda_n/4)=(0,1)$. Let $\mu$
be an independent random cascade. Fix $\kappa>0$. For $j\ge 2$ let
$\mathcal{S}_j=j\log (j)^{-\kappa}$ and $\rho_j=\log (j)^{\alpha}$
with $\alpha>1$. Assume that (\ref{assump}) holds.

For every $q\in J$ and $\xi>1$, with probability one, the property
$\mathcal{P}(\xi,q)$ holds, where $\mathcal{P}(\xi,q)$ is: For
$\mu_q$-almost every $t$, there are infinitely many $k\ge 1$ such that
$\mathcal{B}_{k,1/2} (t)\neq\emptyset$ and there exists
$u\in\{v\in\mathcal{A}^*: \exists\ I\in \mathcal{B}_k^\xi (t),\
I=I_v\}$ such that
\begin{eqnarray}
\label{uu}
GS\big(\mu_q^{u},\mu_q^{u},\tau_\mu'(q)q-\tau_\mu(q), N,\widetilde
\varepsilon\big ) \le \mathcal{S}_{|u|},\quad \mbox{and}\quad
\Vert\mu_q^u\Vert\ge b^{-\rho_{|u|}}.
\end{eqnarray}
\end{theorem}
%%%%%%%%%%%%%%%%%%%%%%%%%%%%%%%%%%%%%%%%%%%%%

%%%%%%%%%%%%%%%%%%%%%%%%%%%%%%%%%%%%%%%%%%%%%
\begin{remark}\label{jo}
The control of $GS\big(\mu_q^{w},\mu^{w},\tau_\mu'(q), N,\widetilde
\varepsilon\big )$ is not useful in deriving (\ref{DIM}).

The result in \cite{MOIBARRALNU} concerning ubiquity conditioned by
Mandelbrot measures invokes a slightly different version of
Theorem~\ref{analogue2}. The proof of this other version is easily
deduced from that of Theorem~\ref{analogue2}.
\end{remark}
%%%%%%%%%%%%%%%%%%%%%%%%%%%%%%%%%%%%%%%%%%%%%
%%%%%%%%%%%%%%%%%%%%%%%%%%%%%%%%%%%%%%%%%%%%%
\begin{proof} 
For $k\ge 1$ and $w\in\mathcal{A}^{k+3}$, notice that
$\mathcal{B}_{k,1/4}(t)\subset \mathcal{B}_{k,1/2}(s)$ for all $t,s\in
I_w$. Let $\mathcal{R}_w=\{n: \exists t\in I_w,\ B(x_n,\lambda_n)\in
\mathcal{B}_{k,1/4}(t)\}$.  Define $n(w)=\inf\big \{n:x_n=\min\{x_m:
m\in \mathcal{R}_w\}\big \}$ if $\mathcal{R}_w\neq \emptyset$ and
$n(w)=0$ otherwise.

If $\xi>1$ and $n(w)>0$, let $u(w)$ be the word encoding the
$b$-adic interval of maximal length included in
$B(x_n,\lambda^\xi_n)$ and whose left end point is minimal. If
$\xi>1$ and $n(w)=0$, let $u(w)$ be the word of generation
$[\xi|w|]$ with prefix $w$ and its $[\xi|w|]-|w|$ last digits
equal to 0. 

%We define so thoroughly $n(w)$ and $u(w)$ because we want them to be
%measurable functions if $\{(x_n,\lambda_n)\}$ is a random set
%containing random variables defined on the same probability space as
%the random variables defining $\mu$.

Now, $w^{(j)}(t)$ being defined as in the statement of
Theorem~\ref{cascades_ren2}, we prove a slightly stronger result than
Theorem~\ref{analogue2}: For every $q\in J$ and $\xi>1$, with
probability one, the property $\widetilde {\mathcal{P}}(\xi,q)$ holds,
where $\widetilde{\mathcal{P}}(\xi,q)$ is: For $\mu_q$-almost every
$t$, if $j$ is large enough, for all $k\ge j$ such that $n\big
(w_{k+3}(t)\big )>0$, $u=u(w_{k+3}(t))$ satisfies (\ref{uu}).

In the sequel we denote $u (w_{k+3}(t) )$ by $u_{k,\xi}(t)$.

\medskip

We fix $\xi>1$ and $q\in K$.

For $j\ge 0$ and $n\ge 1$ define on $\Omega\times [0,1)$ the
random variables
\begin{eqnarray*}
\mathbf{U}^{(j)}(\omega,t)&=&\Vert\mu_q^{u_{j,\xi}(t)}(\omega)\Vert.\\
\mbox{and } \ \
\mathbf{V}^{(j)}_n(\omega,t)&=&S_n^{N,\varepsilon_n,\eta_n} \big
(\mu_q^{u_{j,\xi}(t)}(\omega) , \mu_q^{u_{j,\xi}(t)} (\omega),
q\tau_\mu'(q)-\tau_\mu'(q)\big ).
\end{eqnarray*}
We can use the proof of Proposition~\ref{ren} to deduce that it is
enough to prove
\begin{eqnarray}
\label{33}&& \sum_{j\ge 1}\mathcal{Q}_q\Big (\Big \{\exists \ k\ge
j, \ b^{\rho_{|u_{k,\xi}(t)|}}\sum_{n\ge
\mathcal{S}_{|u_{k,\xi}(t)|}}\mathbf{V}_n^{(k)}(\omega,t)\ge 1/2\Big
\}\Big )<\infty\\
\label{44}\mbox{and}&& \sum_{j\ge 1}\mathcal{Q}_q\Big (\Big \{\exists
\ k\ge j, \ \mathbf{U}^{(j)}(\omega,t) \le
b^{-\rho_{|u_{k,\xi}(t)|}}\Big \}\Big )<\infty.
\end{eqnarray}
Since there exist $c>c'>0$ such that $c'\xi k\le |u_{k,\xi}(t)|\le
c\xi k$ for all $t$, denoting $\bar k=[c\xi k]+1$ and $\widetilde
k=[c'\xi k]$, in order to get (\ref{33}) and (\ref{44}), it is enough
to show that
$$
\begin{cases}
\mathcal{T}=\sum_{j\ge 1}\sum_{k\ge j}\mathcal{Q}_q\big (b^{\rho_{\bar
k}}\sum_{n\ge \mathcal{S}_{\widetilde
k}}\mathbf{V}_n^{(k)}(\omega,t)\ge 1/2\big )<\infty,\\
%%%%%%%%%%
\mathcal{T}'=\sum_{j\ge 1}\sum_{k\ge j}\mathcal{Q}_q(\mathbf{U}^{(k)}\le
b^{-\rho_{\widetilde k}})<\infty.
\end{cases}
$$
Notice that $\displaystyle \mathcal{T}\le 2^h \sum_{j\ge 1}\sum_{k\ge
j}\sum_{n\ge \mathcal{S}_{\widetilde k}}b^{\rho_{\bar
k}h}\mathbb{E}_{\mathcal{Q}_q}\left (\big (\mathbf{V}^{(k)}_n\big
)^h\right )$ if $h\in (0,1)$.

\noindent Mimicking the computations performed in the proof of
Theorem~\ref{cascades_ren2}, one gets
%%%%%%%%%%%%%%%%%%%%%%%%%%%%%%%%%%%%%%%%%%%%%
\begin{eqnarray*}
\int_{[0,1]}\big (V^{u_{k,\xi}(t)}_n(\omega)\big )^h\,
\mu_q(\omega)(dt) =\sum_{w\in\mathcal{A}^{k+3}}\big
(\prod_{k=0}^{k-1}W_{q,w_{k+1}}( w|k)\big ) \big (V_n^{u(w)}\big )^h\,
\Vert\mu_q^{w}\Vert.
\end{eqnarray*}
%%%%%%%%%%%%%%%%%%%%%%%%%%%%%%%%%%%%%%%%%%%%%
Using the independences as well as $p$ and $h$ as in the proof of
Theorem~\ref{cascades_ren2}, one obtains
%%%%%%%%%%%%%%%%%%%%%%%%%%%%%%%%%%%%%%%%%%%%%
$$
\mathbb{E}_{\mathcal{Q}_q}\big(\big (\mathbf{V}^{(k)}_n\big
)^h\big)\le (\mathbb{E}\big (V_n^{u(w)}\big )\big )^{h}\,
\mathbb{E}\big (\Vert\mu_q^w\Vert^p\big )^{{1}/{p}}
$$
where $w$ is any element of $\mathcal{A}^*$. Then our choice for
$\rho_j$ and $\mathcal{S}_j$ ensures that $\mathcal{T}$ is finite.

%%%%%%%%%%%%%%%%%%%%%%%%%%%%%%%%%%%%%%%%%%%%%

\smallskip

Now, for any $h'>0$ one has $\mathcal{T}' \le \sum_{j\ge 1}\sum_{k\ge
j}b^{-\rho_{\widetilde k} h}\mathbb{E}_{\mathcal{Q}_q}\big(\big
(\mathbf{U}^{(k)}\big )^{-h'}\big)$. A computation similar to the
previous one yields, with the same $h$ and $p$, $\displaystyle
\mathbb{E}_{\mathcal{Q}_q}\big(\big (\mathbf{U}^{(k)}\big
)^{-h'}\big)\le \big (\mathbb{E}\big (Y_q^{u(w)}\big )^{-h'/h}\big
)^{h}\big (\Vert\mu_q^w\Vert^p\big )^{{1}/{p}}$ for any element $w$ of
$\mathcal{A}^*$. If $h'$ is chosen small enough, by
Lemma~\ref{supmoments} the right hand side is bounded independently of
$k$ and the conclusion follows from our choice for $\rho_j$.
\end{proof}
%%%%%%%%%%%%%%%%%%%%%%%%%%%%%%%%%%%%%%%%%%%%%
%%%%%%%%%%%%%%%%%%%%%%%%%%%%%%%%%%%%%%%%%%%%%

%%%%%%%%%%%%%%%%%%%%%%%%%%%%%%%%%%%%%%%%%%%%%%%%%%%%%%%%%%%%%%%%%%%
%%%%%%%%%%%%%%%%%%%%%%%%%%%%%%%%%%%%%%%%%%%%%%%%%%%%%%%%%%%%%%%%%%%
%%%%%%%%%%%%%%%%%%%%%%%%%%%%%%%%%%%%%%%%%%%%%%%%%%%%%%%%%%%%%%%%%%%

%%%%%%%%%%%%%%%%%%%%%%%%%%%%%%%%%%%%%%%%%%%%%%%%

%%%%%%%%%%%%%%%%%%%%%%%%%%%%%%%%%%%%%%%%%%%%%%%%
%%%%%%%%%%%%%%%%%%%%%%%%%%%%%%%%%%%%%%%%%%%%%%%%

%%%%%%%%%%%%%%%%%%%%%%%%%%%%%%%%%%%%%%%%%%%%%%%%
%%%%%%%%%%%%%%%%%%%%%%%%%%%%%%%%%%%%%%%%%%%%%%%%
%%%%%%%%%%%%%%%%%%%%%%%%%%%%%%%%%%%%%%%%%%%%%%%%

\end{document}